\documentclass[12pt]{amsart}
\usepackage{mathrsfs}
\usepackage{txfonts}
\usepackage{amssymb}
\makeatletter \@mparswitchfalse \makeatother
\def\eqref#1{(\ref{eq#1})}

\numberwithin{equation}{section}
\begin{document}
\title{Pencils of pairs of projections}
\author{Miaomiao Cui}
  \address{School of Mathematics and Information Science,
  Shaanxi Normal University,
  Xi'an, 710119, People's  Republic of  China}
  \email{cuiye@snnu.edu.cn}
\author{Guoxing Ji$^*$}
  \address{School of Mathematics and Information Science,
  Shaanxi Normal University,
  Xi'an, 710119, People's  Republic of  China}
  \email{gxji@snnu.edu.cn}
\thanks{This research was
supported by the National Natural Science Foundation of China(No. 11771261).}
\thanks{$^*$Corresponding author}
    \maketitle
\begin{abstract}
Let $T$ be a self-adjoint operator on a complex  Hilbert space $\mathcal{H}$. We give a sufficient and necessary condition for $T$ to be the pencil $\lambda P+Q$ of a pair $( P, Q)$ of projections at some point $\lambda\in\mathbb{R}\backslash\{-1, 0\}$. Then we represent all pairs $(P, Q)$ of projections such that $T=\lambda P+Q$ for a fixed $\lambda$, and find that all such pairs are connected if $\lambda\in\mathbb{R}\backslash\{-1, 0, 1\}$. Afterwards, the von Neumann algebra generated by such pairs $(P,Q)$ is characterized. Moreover, we prove that there are at most two real numbers such that $T$ is the pencils at these real numbers for some pairs of projections. Finally, we determine when the real number is unique.\\\\
\subjclass{\textbf{2010 MSC:} 47B15, 47A05, 47C15.}\\
\keywords{\textbf{Keywords:} Pencil, pair of projections,  self-adjoint operator, von Neumann algebra.}
\end{abstract}
\
\baselineskip18pt

\section{Introduction}

The study of pairs of projections has a long history dating back to the late 1940's. In 1948, Dixmier \cite {di} showed that the self-adjoint operator $P+Q$ is a complete unitary invariant of the pair $(P, Q)$ of projections when the pair is in generic position. Subsequently, Davis \cite{da} characterized the self-adjoint operators which are a difference of two orthogonal projections. Afterwards, Kato \cite{ka, ka1} proved that two orthogonal projections are unitarily equivalent when the norm of their difference  is strictly less than 1. Later, Halmos in \cite{ha} represented pairs of projections in generic position by use of graphs of linear transformations. Meanwhile, Fillmore \cite{fi} obtained that any positive invertible operator is a positive linear combination of projections by means of a characterization of operators which are the sum of a pair of projections. On the basis of these remarkable results, several authors recently study the index \cite{am,av}, the difference \cite{an1,ko,sh}, the product \cite{co}, the essential commutators \cite{an2}, some geometric properties(especially minimal geodesics) \cite{an, an3, co1} of pairs of projections and other directions, which leads to questions of unitary equivalence of two pairs of projections.

On the other hand, Raeburn and  Sinclair in \cite{ra} considered when two pairs of projections are unitarily equivalent. They shown that for $\lambda\in(0, +\infty)\backslash\{1\}$, there is a unitary $U$ such that $UPU^{*}=P_{1}$ and $UQU^{*}=Q_{1}$ if and only if the operator $\lambda P+Q$ is unitarily equivalent to the operator $\lambda P_{1}+Q_{1}$. In fact, $\lambda P+ Q$ is said to be the pencil of  $P$ and $Q$ at $\lambda$(\cite[Formula 12.1]{ma}). Among other things, it is natural to consider two problems. One is when a self-adjoint operator $T$ can be expressed as the pencil $\lambda P+Q$ of a pair $(P, Q)$ of projections at some point $\lambda\in\mathbb{R}\backslash\{0\}$. Moreover, a description of all pairs $(P, Q)$ of projections satisfying $T=\lambda P+Q$ for a fixed $\lambda$ is asked for, and the von Neumann algebra generated by those pairs $(P, Q)$ needs to be known. The other is how many real numbers there are such that $T$ is the pencils at these real numbers for pairs of projections, and when the real number is unique. We consider these two problems in this paper.

We now recall some notations. Let $\mathcal{H}$ be an infinite dimensional complex Hilbert space and let $\mathcal{B(H)}$ be the algebra of all bounded linear operators on $\mathcal{H}.$ Denote by $I$ the identity on any Hilbert space and by $\dim\mathcal{M}$ the dimension of  a closed subspace $\mathcal{M}\subseteq\mathcal{H}$. For $T\in \mathcal{B(H)}$, $T^*,~\mathcal{N}(T),~\mathcal{R}(T),$ $\sigma(T)$ and $\sigma_{p}(T)$ stand for the adjoint, the null space, the range, the spectrum and the point spectrum of $T$, respectively. The set of all self-adjoint operators of $\mathcal{B(H)}$ is given by $\mathcal{B(H)}_{s}=\{T\in\mathcal{B(H)}:T^{*}=T\}$ and the set of all orthogonal projections of $\mathcal{B(H)}$ is given by $\mathcal{P(H)}=\{P\in\mathcal{B(H)}_{s}:P=P^2\}.$  For the sake of convenience, by a projection we always mean an orthogonal projection. For $T\in\mathcal{B(H)}_{s}$, $T$ is called a positive operator if $\sigma(T)\subseteq[0, +\infty)$ and denoted by $T\geq0$. Suppose $T_{+}$ and $T_{-}$, respectively, are the positive part and the negative part of self-adjoint operator $T$. For $A,B\in\mathcal B(\mathcal H)_s$, $A\leq B$(resp.  $A<B$) will mean that $ B-A$ is positive (resp. injective positive). For any subset $\mathcal{F}$ of $\mathcal{B(H)}$, denote by $W^{*}(\mathcal F)$ the von Neumann algebra generated by $\mathcal F$. The commutant of $\mathcal{F}$ is given by $\mathcal {F}^{\prime}= \{T\in \mathcal{B(H)}: TF=FT \mbox{ for  all  } F \in\mathcal F\}$. As we all know, $\{T\}^{\prime}$ is a von Neumann algebra for any self-adjoint operator $T$. Let $\mathcal R$ be a von Neumann algebra and let $\mathcal{U(\mathcal R)}$ be the set of all unitary operators in $\mathcal R$. The von Neumann algebra $M_2(\mathcal R)$ and the 2-fold inflation of $\mathcal R$ are given by $M_2(\mathcal R)=\{(R_{ij})_{2\times 2}:R_{ij}\in\mathcal R\}$ and $\mathcal R^{(2)}=\{R\oplus R:R\in\mathcal R\}$. For a subset $S\subseteq\mathbb{R}$, min$S$ and max$S$ stand for the minimum element and the maximum element in $S$, respectively.\\

Before stating the results of the paper, we recall the following notations.\\\\
{\bf Definition 1.1}(see \cite[Formula 12.1]{ma})\quad{\it Suppose $P, Q\in\mathcal{B(H)}$ are two nonzero  projections. We say pencil of pair $(P, Q)$ is  the function $\varphi_{(P, Q)}(\cdot)$ defined as
$$\varphi_{(P, Q)}: \mathbb{R}\longrightarrow \mathcal{B(H)}_{s},~ \varphi(t)=tP+Q.$$}
In particular, we call $\lambda P+Q$ the pencil of pair $(P, Q)$ of projections at $\lambda\in\mathbb{R}$. Fix $\lambda\in\mathbb{R}$, set $\mathcal{T}_{\lambda}=\{(P, Q)\in \mathcal{P(H)}\times\mathcal{ P(H)}: T=\lambda P+Q\}.$\\

Assume that $T$ is a self-adjoint operator and $\lambda\in\mathbb{R}$.  Then   $\mathcal{N}(T)$, $\mathcal{N}(T-I)$, $\mathcal{N}(T-\lambda I)$, $\mathcal{N}(T-(1+\lambda)I)$ and the orthogonal complement $\mathcal{H}_{0}$ of the sum of these reduce $T$. We have that according to the decomposition
$$\mathcal{H}=\mathcal{N}(T)\oplus \mathcal{N}(T-I)\oplus \mathcal{ N}(T-\lambda I)\oplus \mathcal{N}(T-(1+\lambda)I)\oplus \mathcal{H}_{0},\eqno{(1.1)}$$
$T$ is given by $$T=0\oplus I\oplus \lambda I\oplus(1+\lambda)I\oplus T_0,\eqno{(1.2)}$$
where $T_{0}= T|_{\mathcal{H}_0}$. We call $T_{0}$ the generic part of $T$ with respect to $\lambda$.
 %It is easy to check $\frac{1+|\lambda|}{2}, \frac{|1-|\lambda||}{2}\notin\sigma_{p}(T_{0}-\frac{1+\lambda}{2}I)$.
If
$$\mathcal{N}(T)= \mathcal{N}(T-I)=\mathcal{N}(T-\lambda I)=\mathcal{N}(T-(1+\lambda)I)= \{0\},$$
then we say $T$ is in generic position with respect to $\lambda$.

If $T$ is the pencil of pair $(P, Q)$ of projections at $\lambda$, then (cf. \cite{ha})
$$\left\{
  \begin{array}{ll}
\mathcal{N}(T)=\mathcal{N}(P)\cap \mathcal{N}(Q), \quad\quad \quad ~ \mathcal{N}(T-I)=\mathcal{N}(P)\cap \mathcal{R}(Q), \\
\mathcal{N}(T-\lambda I)=\mathcal{R}(P)\cap \mathcal{N}(Q), \quad \mathcal{N}(T-(1+\lambda)I)=\mathcal{R}(P)\cap \mathcal{R}(Q).
  \end{array}
\right.$$
In this case, it is known that  $T$ is in generic position with respect to $\lambda$ if and only if the pair $(P,Q)$ of projections is in generic position in the sense of \cite{ha}.

The contents of the paper are as follows. In Section 2, we establish a characterization for $T$ to be the pencil of a pair of projections at some point $\lambda\in\mathbb{R}\backslash\{-1, 0\}$. By use of this  characterization, one can obtain the general representations of all pairs $(P, Q)$ of projections such that $T=\lambda P+Q$ for a fixed $\lambda$, and find that all such pairs are connected if $\lambda\in\mathbb{R}\backslash\{-1, 0, 1\}$. Furthermore, we characterize the von Neumann algebra generated by all pairs $(P, Q)\in\mathcal{ T}_{\lambda}$ of projections as well as its commutant. Section 3 is devoted to show that there are at most two real numbers such that $T$ is the pencils at these real numbers for some pairs of projections. In addition, we determine when the real number is unique.
\section{When T is the pencil of a pair of projections at some point}
Let $T\in\mathcal{ B(H)}$ be a self-adjoint operator.  We have known  when $T$ is a difference of two projections(cf. \cite{da,sh}).  We consider when $T$ is the pencil of two projections at some point $\lambda \in\mathbb R\backslash\{-1, 0\}$. Put $\Lambda=\mathbb{R}\backslash\{-1, 0\}$ and $\Lambda_1=\mathbb{R}\backslash\{-1, 0, 1\}$. Suppose $T$ is in generic position with respect to $\lambda\in\Lambda$. Next, we consider when $T$ is the pencil of a pair of projections at $\lambda$, and then describe all such pairs. Afterwards, we extend the results to the general case.\\\\
{\bf Proposition 2.1}\quad{\it Let $T\in\mathcal{ B(H)}$ be a self-adjoint operator in generic position with respect to $\lambda\in\Lambda$. Then $T$ is the pencil of a pair of projections at $\lambda$ if and only if there is a positive  operator $B$ on a Hilbert space $\mathcal K$ satisfying $\frac{|1-|\lambda||}{2}I<B<\frac{1+|\lambda|}{2}I$ such that $T-\frac{1+\lambda}{2}I$ is unitarily equivalent to
 $$B\oplus (-B). \eqno{(2.1)}$$}
{\bf Proof.}\quad Suppose $T=\lambda P+Q$, where $(P,~Q)$ is a pair of projections. It is elementary   that $P$ and $Q$ are in generic position(\cite{ha}). Firstly, we conclude $$\frac{(1-|\lambda|)^{2}}{4}I<(T-\frac{1+\lambda}{2}I)^{2}<\frac{(1+|\lambda|)^{2}}{4}I.\eqno{(2.2)}$$
Indeed, straightforward calculation shows that $(2.2)$ is equivalent to
$$-
\frac{|\lambda|+\lambda}{2}I<-\lambda(P-Q)^{2}<\frac{|\lambda|-
\lambda}{2}I,$$
which is also equivalent to $0<(P-Q)^{2}<I.$
From \cite[Proposition 4]{am}, we know $0<(P-Q)^{2}\leq\|P-Q\|^{2}I<I$ since $P$ and $Q$ are in generic position. Hence $\sigma(T-\frac{1+\lambda}{2}I)\subseteq [-\frac{1+|\lambda|}{2},-\frac{|1-|\lambda||}{2}]\cup [\frac{|1-|\lambda||}{2},\frac{1+|\lambda|}{2}]$. Let $T-\frac{1+\lambda}{2}I=\int \lambda dE_{\lambda}$ be the spectral decomposition of $T-\frac{1+\lambda}{2}I$ and let  $\mathcal{H}_{+}=E[\frac{|1-|\lambda||}{2},\frac{1+|\lambda|}{2}]$ and $\mathcal{H}_{-}=E[-\frac{1+|\lambda|}{2},-\frac{|1-|\lambda||}{2}]$.
Apparently,   $T-\frac{1+\lambda}{2}I=(T-\frac{1+\lambda}{2}I)_{+}\oplus(T-\frac{1+\lambda}{2}I)_{-}$ with respect to
$\mathcal{H}=\mathcal{H}_{+}\oplus \mathcal{H}_{-}$.

Set $S=PQ-QP$, it is known that $S$ is injective with dense range. By  a direct computation, we have $(T-\frac{1+\lambda}{2}I)S=-S(T-\frac{1+\lambda}{2}I)$. Thus there exists a unitary operator $U$ such that $U(T-\frac{1+\lambda}{2}I)=-(T-\frac{1+\lambda}{2}I)U$ according to \cite[Corollary 6.11]{con}. Clearly, $U(T-\frac{1+\lambda}{2}I)^{2}=(T-\frac{1+\lambda}{2}I)^{2}U.$ Assume that $U=\left(\begin{array}{cccc}U_{11}&U_{12}\\U_{21}&U_{22}\\\end{array}\right)$ under $\mathcal{H}=\mathcal{H}_{+}\oplus \mathcal{H}_{-}$. From $(2.2)$, we see that  $(T-\frac{1+\lambda}{2}I)_{+}$ and $(T-\frac{1+\lambda}{2}I)_{-}$ are injective. One obtains $U_{11}=U_{22}$ and $U_{21}(T-\frac{1+\lambda}{2}I)_{+}=(T-\frac{1+\lambda}{2}I)_{-}U_{21}$ by calculation, which means that $U_{21}$ is a unitary operator from $\mathcal{H}_{+}$ onto $\mathcal{H}_{-}$. Hence $U_{21}(T-\frac{1+\lambda}{2}I)_{+}U_{21}^{*}=(T-\frac{1+\lambda}{2}I)_{-}$. This implies that $T-\frac{1+\lambda}{2}I$ is unitarily equivalent to $B\oplus-B$, where $B=(T-\frac{1+\lambda}{2}I)_{+}$. It clearly forces $\frac{|1-|\lambda||}{2}I<B<\frac{1+|\lambda|}{2}I$ from $(2.2)$.

Conversely, assume that there exists a positive operator $B$ satisfying $\frac{|1-|\lambda||}{2}I<B<\frac{1+|\lambda|}{2}I$ such that $T-\frac{1+\lambda}{2}I$ is unitarily equivalent to $B\oplus -B$. Without loss of generality, we  may assume that  $T-\frac{1+\lambda}{2}I=B\oplus -B$. Put
$$P=\frac{1}{2\lambda}\left(\begin{array}{cccc}P_{11}&P_{12}\\
P_{12}^{*}&P_{22}\\\end{array}\right)~\hbox{ and } ~Q=\frac{1}{2}\left(\begin{array}{cccc}Q_{11}&Q_{12}\\
 Q_{12}^{*}&Q_{22}\\\end{array}\right),\eqno{(2.3)}$$
where if $\lambda\in\Lambda_{1}$, set
$$\left\{
  \begin{array}{ll}
P_{11}=B^{-1}(B+\frac{1+\lambda}{2}I)(B-\frac{1-\lambda}{2}I),~P_{22}=-B^{-1}(B-\frac{1+\lambda}{2}I)(B+\frac{1-\lambda}{2}I),\\
Q_{11}=B^{-1}(B+\frac{1+\lambda}{2}I)(B+\frac{1-\lambda}{2}I),~Q_{22}=-B^{-1}(B-\frac{1+\lambda}{2}I)(B-\frac{1-\lambda}{2}I),\\
P_{12}=B^{-1}\{-(B^{2}-\frac{(1+\lambda)^{2}}{4}I)(B^{2}-\frac{(1-\lambda)^{2}}{4}I)\}^{\frac{1}{2}},~Q_{12}=-P_{12}
  \end{array}
\right.$$
and if $\lambda=1$, set
$$\left\{
  \begin{array}{ll}
P_{11}=B+\frac{1+\lambda}{2}I,~P_{22}=-(B-\frac{1+\lambda}{2}I),\\
Q_{11}=B+\frac{1+\lambda}{2}I,~Q_{22}=-(B-\frac{1+\lambda}{2}I),\\
P_{12}=\{-(B+\frac{1+\lambda}{2}I)(B-\frac{1+\lambda}{2}I)\}^{\frac{1}{2}}B,~Q_{12}=-P_{12}.
  \end{array}
\right.$$
It is easy to check that $(P, Q)$ is a pair of projections satisfying
\begin{center}
$\lambda P+Q-\frac{1+\lambda}{2}I=B\oplus -B,$
\end{center}
which gives that $T=\lambda P+Q$.\qed\\\\
{\bf Remark 2.2}\quad We may assume that $T-\frac{1+\lambda}{2}I$ has the form $(2.1)$ with $\lambda\in\Lambda$ by Proposition 2.1 if $T$ is in generic position at $\lambda$. That is, $T-\frac{1+\lambda}{2}I=B\oplus(-B),$ where $\frac{|1-|\lambda||}{2}I<B<\frac{1+|\lambda|}{2}I$.
Now for any $U\in\mathcal U(\{B\}^{\prime})$, put
$$P_U=\frac{1}{2\lambda}\left(\begin{array}{cccc}P_{11}&P_{12}U\\
 U^{*}P_{12}&P_{22}\\\end{array}\right) \hbox{ and }
Q_U=\frac{1}{2}\left(\begin{array}{cccc}Q_{11}&Q_{12}U \\
 U^{*}Q_{12}&Q_{22}\\\end{array}\right),\eqno{(2.4)}$$
where $P_{ij}$ and $Q_{ij}(i, j=1, 2)$ are defined as in $(2.3)$. It is known that $T$ is the pencils of pairs $(P_U, Q_U)$ of projections at $\lambda$.\\

Therefore we have the following result.\\\\
{\bf Theorem 2.3}\quad{\it Let $T\in\mathcal{ B(H)}$ be a self-adjoint operator in generic position with respect to $\lambda\in\Lambda$ such that $T-\frac{1+\lambda}{2}I$ has the form $(2.1)$. Then $T$ is the pencil of a pair of projections at $\lambda$. Moreover, $$\mathcal{T}_\lambda=\{(P_U,~Q_U): U\in\mathcal U(\{B\}^{\prime})\}.$$}\\
{\bf Proof.}\quad
 %Suppose $T-\frac{1+\lambda}{2}=B\oplus (-B)$ be the standard form. Put $P_U$ and $Q_U$ be defined as in $(6)$ and $(7)$. Then $T=\lambda P_U+Q_U$.
Let $\frac{|1-|\lambda||}{2}I<B<\frac{1+|\lambda|}{2}I$ such that $T-\frac{1+\lambda}{2}I=B\oplus (-B)$ in terms of $\mathcal{H}=\mathcal{K}\oplus\mathcal{K}$. Suppose $T$ is the pencil of pair $(P, Q)$ of projections at $\lambda$ and
\begin{center}$P=\left(\begin{array}{cc}P_{11}&P_{12}\\P_{12}^*&P_{22}
\end{array}\right),~~
Q=\left(\begin{array}{cc}Q_{11}&Q_{12}\\Q_{12}^*&Q_{22}
\end{array}\right).$
\end{center}
with respect to $\mathcal{H}=\mathcal{K}\oplus\mathcal{K}$. Then
\begin{center}
$T=\left(\begin{array}{cccc}B&0\\0&-B\\\end{array}\right)+\frac{1+\lambda}{2}I= \left(\begin{array}{cc}\lambda P_{11}+Q_{11}&\lambda P_{12}+Q_{12}\\ \lambda P_{12}^*+Q_{12}^* & \lambda P_{22}+Q_{22}
\end{array}\right).$
\end{center}
Clearly,
$$\left\{\begin{array}{l}\lambda P_{11}+Q_{11}=B+\frac{1+\lambda}{2}I,\\
\lambda P_{12}+Q_{12}=0,\\ \lambda P_{22}+Q_{22}=-B+\frac{1+\lambda}{2}I.\end{array}\right.\eqno{(2.5)}$$
Note that $P[T^{2}-(1+\lambda)T]=\lambda(PQP-P)=[T^{2}-(1+\lambda)T]P$ and
$$T^{2}-(1+\lambda)T+\frac{(1+\lambda)^{2}}{4}I=(T-\frac{1+\lambda}{2}I)^{2}=B^{2}\oplus B^{2}.  $$  It follows that $P|T-\frac{1+\lambda}{2}I|=|T-\frac{1+\lambda}{2}I|P$. It is immediate that
$$\left(\begin{array}{cc}P_{11}B&P_{12}B\\P_{12}^*B & P_{22}B\end{array}\right)=\left(\begin{array}{cccc}BP_{11}&BP_{12}\\
 BP_{12}^* & BP_{22} \end{array}\right),$$
which implies that $P_{ij}\in\{B\}^{\prime}$ and $Q_{ij}\in\{B\}^{\prime}$ by $(2.5)$ for $i,j=1,2$.
Since $P$ and $Q$ are projections, we see that $$\left\{\begin{array}{l}P_{11}^{2}+P_{12}P_{12}^*= P_{11},~~~~
P_{12}^*P_{12}+P_{22}^{2}= P_{22},\\Q_{11}^{2}+Q_{12}Q_{12}^*=Q_{11},~~Q_{12}^*Q_{12}+Q_{22}^{2}= Q_{22}.\end{array}\right.\eqno{(2.6)}$$
From $(2.5)$ and $(2.6),$  we then obtain
$$\left\{\begin{array}{l}(B+\frac{1+\lambda}{2}I-\lambda P_{11})^{2}+\lambda^{2}P_{12}P_{12}^{*}=B+\frac{1+\lambda}{2}I-\lambda P_{11},\\(-B+\frac{1+\lambda}{2}I-\lambda P_{22})^{2}+\lambda^{2}P_{12}^{*}P_{12}=-B+\frac{1+\lambda}{2}I-\lambda P_{22},\end{array}\right.$$
which yields $$\left\{\begin{array}{l}(B+\frac{1+\lambda}{2}I)(B-\frac{1-\lambda}{2}I)=2\lambda BP_{11},\\(B-\frac{1+\lambda}{2}I)(B+\frac{1-\lambda}{2}I)=-2\lambda BP_{22}.\end{array}\right.\eqno{(2.7)}$$
In fact, $B$ is injective if $\lambda=1$, and $B$ is invertible if $\lambda\in\Lambda_{1}$. We thus get
\begin{center}
$P_{11}=\frac{1}{2}(I+B),~ P_{22}=\frac{1}{2}(I-B)$
\end{center}
when $\lambda=1$ and
$$\left\{\begin{array}{l} P_{11}=\frac{1}{2\lambda}B^{-1}(B+\frac{1+\lambda}{2}I)(B-\frac{1-\lambda}{2}I),\\
P_{22}=-\frac{1}{2\lambda}B^{-1}(B-\frac{1+\lambda}{2}I)(B+\frac{1-\lambda}{2}I)
\end{array}\right.$$
when $\lambda\in\Lambda_{1}$ from $(2.7)$. According to $(2.5)$ again, it follows that
\begin{center}
$Q_{11}=\frac{1}{2}(I+B),~Q_{22}=\frac{1}{2}(I-B)$
\end{center}
when $\lambda=1$ and
$$\left\{\begin{array}{l} Q_{11}=\frac{1}{2}B^{-1}(B+\frac{1+\lambda}{2}I)(B+\frac{1-\lambda}{2}I),\\
Q_{22}=-\frac{1}{2}B^{-1}(B-\frac{1+\lambda}{2}I)(B-\frac{1-\lambda}{2}I)
\end{array}\right.$$
when $\lambda\in\Lambda_{1}$.
These imply that
$$\left\{\begin{array}{l} P_{12}P_{12}^*= P_{11}(I-P_{11})= \frac{1}{4}(I-B^{2})=P_{22}(I-P_{22})=P_{12}^*P_{12}\quad\hbox{ if }\lambda=1,\\
P_{12}P_{12}^*= P_{11}(I-P_{11})= \frac{1}{4\lambda^{2}}B^{-2}[-(B^{2}-\frac{(1+\lambda)^{2}}{4}I)(B^{2}-\frac{(1-\lambda)^{2}}{4}I)]\\

\quad\quad\quad\quad\quad\quad\quad\quad\quad\quad  =P_{22}(I-P_{22})=P_{12}^*P_{12}\quad\hbox{ if }\lambda\in\Lambda_{1} \end{array}\right.$$
by $(2.6)$. The fact that $\frac{|1-|\lambda||}{2}I<B<\frac{1+|\lambda|}{2}I$ ensures that $P_{12}$ is injective normal operator and
$$\left\{\begin{array}{l} |P_{12}|=\frac{1}{2}(I-B^{2})^{\frac{1}{2}}\quad\hbox{ if }\lambda=1,\\
|P_{12}|=\frac{1}{2\lambda}B^{-1}[-(B^{2}-\frac{(1+\lambda)^{2}}{4}I)(B^{2}-\frac{(1-\lambda)^{2}}{4}I)]^{\frac{1}{2}}\quad\hbox{ if }\lambda\in\Lambda_{1}.
\end{array}\right.$$
Let $P_{12}=U|P_{12}|$ be the polar decomposition of $P_{12}$, where $U\in\mathcal{ B(K)}$ is a unitary  operator  in $\{|P_{12}|\}^{\prime}$. Combining these with $(2.5)$, we have
$$\left\{\begin{array}{l}Q_{12}=-\lambda P_{12}=-\frac{1}{2} U[-(B^{2}-\frac{(1+\lambda)^{2}}{4}I)]^{\frac{1}{2}}\quad\hbox{ if }\lambda=1,\\
Q_{12}=-\lambda P_{12}=-\frac{1}{2} UB^{-1}[-(B^{2}-\frac{(1+\lambda)^{2}}{4}I)(B^{2}-\frac{(1-\lambda)^{2}}{4})I]^{\frac{1}{2}}\quad \hbox{ if }\lambda\in\Lambda_{1}.
\end{array}\right.$$

Next, we only need to prove $U\in\{B\}^{\prime}$, which finishes the proof of the theorem. Indeed, if $\lambda=1$, it is clear that $U\in\{B\}^{\prime}$. If $\lambda\in\Lambda_{1}$, it follows that
$UP_{12}^*P_{12}=P_{12}^*P_{12}U$ from $U\in\{|P_{12}|\}^{\prime}$. By an elementary calculation, one knows that
%$$ UB^{-2}(B^{2}-\frac{(1+\lambda)^{2}}{4}I)(B^{2}-\frac{(1-\lambda)^{2}}{4}I)=B^{-2}(B^{2}-\frac{(1+\lambda)^{2}}{4}I)(B^{2}-\frac{(1-\lambda)^{2}}{4}I)U,$$
%which is equivalent to $$U(B-\frac{(1+\lambda)^{2}}{4}B^{-1})(B-\frac{(1-\lambda)^{2}}{4}B^{-1})=(B-\frac{(1+\lambda)^{2}}{4}B^{-1})(B-\frac{(1-\lambda)^{2}}{4}B^{-1})U.$$
$U(B^{2}+\frac{(1-\lambda^{2})^{2}}{4}B^{-2})=(B^{2}+\frac{(1-\lambda^{2})^{2}}{4}B^{-2})U,$ which gives
\begin{center}
$U(B+\frac{|1-\lambda^{2}|}{2}B^{-1})^{2}=(B+\frac{|1-\lambda^{2}|}{2}B^{-1})^{2}U.$
\end{center}
As $B+\frac{|1-\lambda^{2}|}{2}B^{-1}=B^{-1}(B^{2}+\frac{|1-\lambda^{2}|}{2}I)$ is an invertible positive operator, we have $U(B+\frac{|1-\lambda^{2}|}{2}B^{-1})=(B+\frac{|1-\lambda^{2}|}{2}B^{-1})U.$
Set
\begin{center}
$f(t)=(\frac{t^{2}}{4}-\frac{|1-\lambda^{2}|}{2})^{\frac{1}{2}}+\frac{t}{2},~~ t\in\sigma(B+\frac{|1-\lambda^{2}|}{2}B^{-1}).$
\end{center}%t\in[-\frac{1+|\lambda|}{2},-\frac{|1-|\lambda||}{2}]\cup [\frac{|1-|\lambda||}{2},\frac{1+|\lambda|}{2}].$$
By functional calculus, $B=f(B+\frac{|1-\lambda^{2}|}{2}B^{-1})$, it induces that $U\in\{B\}^{\prime}$.\qed \\

Fix $\lambda\in\Lambda$, Theorem 2.3 states that there is a one-to-one correspondence between $\mathcal{T}_{\lambda}$ and $\mathcal U(\{B\}^{\prime})$ if $T-\frac{1+\lambda}{2}I$ has the form $(2.1)$. Consequently, one can obtain the following result.\\\\
{\bf Corollary 2.4}\quad{\it Let $T\in\mathcal{ B(H)}$ be a self-adjoint operator in generic position with respect to $\lambda\in\Lambda$. If $T-\frac{1+\lambda}{2}I$ has the form $(2.1)$, then $\mathcal{T}_\lambda$ is a connected subset of $\mathcal{P(H)}\times\mathcal{ P(H)}$.}\\

Assume that $T\in\mathcal{ B(H)}$ is a self-adjoint operator such that $T-\frac{1+\lambda}{2}I$ has the form $(2.1)$. Obviously, $P_U$ and $Q_U$ do not commute from Theorem 2.3. It is known that a von Neumann algebra $\mathcal{R}$ is generated by all projections in $\mathcal{R}$. Next, we consider the von Neumann algebra $W^*(\mathcal{T}_\lambda)$ generated by all pairs of projections in $\mathcal{T}_\lambda$.\\\\
{\bf Theorem 2.5}\quad{\it Let $T\in\mathcal{ B(H)}$ be a self-adjoint operator in generic position with respect to $\lambda\in\Lambda$. If $T-\frac{1+\lambda}{2}I$ has the form $(2.1)$, then $$W^*(\mathcal{T}_\lambda)=M_{2}(\{B\}^{\prime}) \mbox{ and } (W^*(\mathcal{T}_\lambda))^{\prime}= {W^*(B})^{(2)}.$$}
{\bf Proof.}\quad Assume that $\frac{|1-|\lambda||}{2}I<B<\frac{1+|\lambda|}{2}I$ such that $T-\frac{1+\lambda}{2}I=B\oplus (-B)$ in terms of $\mathcal{H}=\mathcal{K}\oplus\mathcal{K}$. Let us firstly prove $(W^*(\mathcal{T}_\lambda))^{\prime}=W^*(B)^{(2)}$. Suppose $A\in (W^*(\mathcal{T}_\lambda))^{\prime}$ and $A=\left(\begin{array}{cccc}A_{11}&A_{12}\\
 A_{21}&A_{22}\\\end{array}\right)$ with respect to  $\mathcal{H}=\mathcal{K}\oplus\mathcal{K}$. This gives $A(T-\frac{1+\lambda}{2}I)=(T-\frac{1+\lambda}{2}I)A$, which yields $A(T-\frac{1+\lambda}{2}I)^2=(T-\frac{1+\lambda}{2}I)^2A.$ From that $B$ is injective with dense range,
it induces that $A_{ij}=0$ for $i,j=1,2$ and $A_{ii}B=BA_{ii}$ for $i\not=j$ by direct calculation.
Note that $AP_{U}=P_{U}A$ for every unitary operator $U\in\{B\}^{\prime}$, we easily obtain
$$\left(\begin{array}{cccc}A_{11}P_{11}&A_{11}UP_{12}\\
A_{11}U^{*}P_{12}&A_{22}P_{22}\\\end{array}\right)=\left(\begin{array}{cccc}P_{11}A_{11}&UP_{12}A_{22}\\
U^{*}P_{12}A_{11}&P_{22}A_{22}\\\end{array}\right),$$
and so $A_{11}UP_{12}=UP_{12}A_{22}.$ It follows from the fact that $P_{12}$ is injective that
$A_{11}U=UA_{22}$ for all unitary operators $U\in\{B\}^{\prime}$. This clearly forces $A_{11}=A_{22}$ and $A_{11}U=UA_{11}$ for all unitary operators
$U\in \{B\}^{\prime}$. We see at once that $A_{11}\in W^*(B)$, and consequently $(W^*(\mathcal{T}_\lambda))^{\prime} =(W^*(B))^{(2)}$.
Using the double commutant theorem, $W^*(\mathcal{T}_\lambda)=W^*(\mathcal{T}_\lambda)^{\prime\prime}=M_{2}( \{B\}^{\prime}).$\qed\\

Let $T\in\mathcal{ B(H)}$ be a self-adjoint operator in generic position with respect to $\lambda\in\Lambda$. If $T$ is the pencil of a pair of projections at $\lambda$, it follows from Theorem 2.5 and Proposition 2.1 that $W^*(\mathcal{T}_\lambda)$ and $(W^*(\mathcal{T}_\lambda))^{\prime}$ are unitarily equivalent to
$M_{2}(\{(T-\frac{1+\lambda}{2}I)_+\}^{\prime})$ and ${W^*((T-\frac{1+\lambda}{2}I)_+})^{(2)}$, respectively.\\

Suppose $T\in\mathcal{B(H)}$ is a self-adjoint operator having a generic part $T_{0}$ with respect to $\lambda\in\Lambda$. Then $T=0\oplus I\oplus \lambda I\oplus (1+\lambda)I\oplus T_0$ in terms of $(1.1)$ from $(1.2)$. It is immediate that $T$ is the pencil of a pair of projections at $\lambda$ if and only if so is $T_0$. Next, we shall give the representations of all pairs $(P, Q)$ of projections such that $T$ is the pencils of pairs $(P, Q)$ of projections at a fixed $\lambda\in\Lambda$. Without loss of generality, suppose $T_0-\frac{1+\lambda}{2}I$ has the form $(2.1)$. That is, $T_0-\frac{1+\lambda}{2}I=B\oplus (-B)$
with $\frac{|1-|\lambda||}{2}I<B<\frac{1+|\lambda|}{2}I$. For any unitary operator $U\in \mathcal U(\{B\}^{\prime})$ and any projection $E\in \mathcal{P(N}(T-I))$, we define
$$\left\{\begin{array}{l}
P_{\lambda, U}=0\oplus  E\oplus I\oplus P_U,~Q_{\lambda, U}=0\oplus(I-E)\oplus I\oplus Q_U\hbox{ if } \lambda=1,\\
P_{\lambda, U}=0\oplus0\oplus  I\oplus I\oplus P_U,~Q_{\lambda, U}=0\oplus I\oplus 0\oplus I\oplus Q_U\hbox{ if } \lambda\in\Lambda_{1}
\end{array}\right.\eqno{(2.8)}$$
with respect to $(1.2)$, where $P_U$ and $Q_U$ are defined as in $(2.4)$.\\\\
{\bf Theorem 2.6}\quad{\it Let $T\in\mathcal{ B(H)}$ be a self-adjoint operator having a generic part $T_0$ with respect to $\lambda\in\Lambda$. If $T_{0}-\frac{1+\lambda}{2}I$ has the form $(2.1)$, then $T$ is the pencil of a pair of projections at $\lambda$. Moreover,
$$\left\{\begin{array}{l}
  \mathcal{T}_\lambda=\{(P_{\lambda, U},Q_{\lambda, U}): E\in \mathcal{P(N}(T-I)), U\in\mathcal U(\{B\}^{\prime})\}\quad \hbox{ if }\lambda=1,\\
 \mathcal{T}_\lambda=\{(P_{\lambda, U},Q_{\lambda, U}): U\in\mathcal{U}(\{B\}^{\prime})\}\quad\hbox{ if } \lambda\in\Lambda_1.\\
\end{array}\right.$$}
{\bf Proof.}\quad Obviously, $T$ is the pencils of pairs $(P_{\lambda, U},Q_{\lambda, U})$ of projections at $\lambda$, where $(P_{\lambda, U},Q_{\lambda, U})$ is defined as in $(2.8)$.

Conversely, assume that $T$ is the pencil of pair $(P,~Q)$ of projections at $\lambda$.
Thus $T_{0}=\lambda P_{0}+Q_{0}$, where $P_{0}=P|_{\mathcal{H}_0}$ and $Q_{0}=Q|_{\mathcal{H}_0}$ are projections. It follows that there is a unitary operator $U\in\{B\}^{\prime}$ such that $P_{0}=P_U$ and $Q_{0}=Q_U$ from Theorem 2.3.

If $\lambda=1$, we would have
\begin{center}
$\left\{\begin{array}{l} P|_{\mathcal{N}(T)}=0=Q|_{\mathcal{N}(T)},~P|_{\mathcal{N}(T-I)}=E,\\
Q|_{\mathcal{N}(T-I)}=I-E,~P|_{\mathcal{N}(T-2I)}=I=Q|_{\mathcal{N}(T-2I)},\end{array}\right.$
\end{center}
which yields $P=0\oplus E\oplus I\oplus P_U=P_{\lambda, U}$ and $Q=0\oplus(I-E)\oplus I\oplus Q_U=Q_{\lambda, U}$.

If $\lambda\in\Lambda_{1}$, then
\begin{center}
$\left\{\begin{array}{l} P|_{\mathcal{N}(T)}=0=Q|_{\mathcal{N}(T)}, P|_{\mathcal{N}(T-I)}=0=Q|_{\mathcal{N}(T-\lambda I)},\\
P|_{\mathcal{N}(T-\lambda I)}=I=Q|_{\mathcal{N}(T-I)}, P|_{\mathcal{N}(T-(1+\lambda)I)}=I=Q|_{\mathcal{N}(T-(1+\lambda)I)}. \end{array}\right.$
\end{center}
Therefore $P=0\oplus0\oplus  I\oplus I\oplus P_U=P_{\lambda, U}$ and $Q=0\oplus I\oplus 0\oplus I\oplus Q_U=Q_{\lambda, U}.$\qed\\

Theorem 2.6 tells us the relationship between the pair $(P, Q)$ and the operator $T=\lambda P+Q$ for $\lambda\in\Lambda$. For any unitary $U\in\{T\}^{\prime}$ and for any pair $(P, Q)\in\mathcal{T}_{\lambda}$, we easily see $(UPU^{*}, UQU^{*})\in\mathcal{T}_{\lambda}$. On the contrary, if $\lambda\in\Lambda_{1}$, then it is worth observing that for any pairs $(P, Q)$, $(P_{1} Q_{1})\in\mathcal{T}_{\lambda}$, there exists a unitary operator $U\in\{T\}^{\prime}$ such that $UPU^{*}=P_{1}$ and $UQU^{*}=Q_{1}$. Consequently, if $\lambda\in\Lambda_{1}$ and  $T=\lambda P+Q$, then
$$\mathcal{T}_{\lambda}=\{(UPU^{*}, UQU^{*}): U\in\{T\}^{\prime}\}.$$

Raeburn and  Sinclair  proved  for $\lambda\in(0, +\infty)\backslash\{1\}$, there is a unitary $U$ on $\mathcal{H}$ such that $UPU^{*}=P_{1}$ and $UQU^{*}=Q_{1}$ if and only if the operator $\lambda P+Q$ is unitarily equivalent to the operator $\lambda P_{1}+Q_{1}$(cf. \cite[Theorem 3.1]{ra}). From Theorem 2.6, we have the following result which extends the case $\lambda\in(0, +\infty)\backslash\{1\}$ to the case $\lambda\in\Lambda_{1}$ by direct calculation.\\\\
{\bf Corollary 2.7}\quad{\it Let $(P, Q)$ and $(P_{1}, Q_{1})$ be two pairs of projections on Hilbert space $\mathcal{H}$ and $\lambda\in\Lambda_{1}$. Then there is a unitary operator $U$ such that $UPU^{*}=P_{1}$ and $UQU^{*}=Q_{1}$ if and only if  $\lambda P+Q$  and $\lambda P_{1}+Q_{1}$ are  unitarily equivalent.}\\

In fact, Corollary 2.7 reveals that the unitary equivalence class of $\lambda P+Q$ uniquely determines that of the pair $(P, Q)$. From $(2.1)$, we furthermore know that the unitary equivalence class of the positive operator $B$ defined as in $(2.1)$ identifies that of the pair $(P, Q)$. We next extend Theorem 2.5 to the general case.\\\\
{\bf Theorem 2.8}\quad{\it  Let $T\in\mathcal{ B(H)}$ be a self-adjoint operator having a generic part $T_0$ with respect to $\lambda\in\Lambda$. If $T_{0}-\frac{1+\lambda}{2}I$ has the form $(2.1)$, then
$$\left\{\begin{array}{l} W^*(\mathcal{T}_\lambda)=\mathbb C I\oplus \mathcal{B(N}(T-I))\oplus \mathbb C I\oplus M_{2}(\{B\}^{\prime}),\\
(W^*(\mathcal{T}_\lambda))^{\prime}=\mathcal{B(N}(T))\oplus\mathbb C I \oplus \mathcal{B(N}(T-2I))\oplus {W^*(B})^{(2)}
\end{array}\right.$$
when $\lambda=1$ and
$$\left\{\begin{array}{l}W^*(\mathcal{T}_\lambda)=\mathbb C I\oplus \mathbb C I\oplus \mathbb C I\oplus \mathbb C I\oplus M_{2}(\{B\}^{\prime}),\\
(W^*(\mathcal{T}_\lambda))^{\prime}=\mathcal{B(N}(T))\oplus\mathcal{B(N}(T-I))\oplus\mathcal{B(N}(T-\lambda I))\\
\quad\quad\quad\quad\quad~\oplus\mathcal{B(N}(T-(1+\lambda)I))\oplus {W^*(B})^{(2)}\end{array}\right.$$ when $\lambda\in\Lambda_{1}$.}\\\\
{\bf Proof.}\quad Assume that $F\in (W^*(\mathcal{T}_{\lambda}))^{\prime}$ is a projection. This gives $FT=TF$. It follows that $\mathcal{N}(T),$ $\mathcal{N}(T-I),$ $\mathcal{N}(T-\lambda I)$, $\mathcal{N}(T-(1+\lambda)I)$ and $\mathcal{H}_0$ reduce $F$. Suppose $F=F_1\oplus F_2\oplus F_3\oplus F_4\oplus F_5$ for some projections $F_i$, $i=1,...,5$. Note that $FP_{\lambda,U}=P_{\lambda,U}F$ and $ FQ_{\lambda,U}=Q_{\lambda,U}F$ for all $(P_{\lambda,U},Q_{\lambda,U})\in \mathcal{T}_{\lambda}$ by Theorem 2.6.

If $\lambda=1$, then $(F_{2}\oplus F_{3})E= E(F_{2}\oplus F_{3})$, $F_{5}P_U= P_UF_{5}$ and $F_{5}Q_U= Q_UF_{5}$ for all $E\in \mathcal{P(N}(T-I))$ and $U\in\mathcal U(\{B\}^{\prime})$. This implies that $F_2\oplus F_{3}=0$ or $I$ and $F_5\in (W^*(B))^{(2)}$ from Theorem 2.5. If $\lambda\in\Lambda_{1}$, we have $F_{5}P_U= P_UF_{5}$ and $F_{5}Q_U= Q_UF_{5}$ for
all $U\in\mathcal U(\{B\}^{\prime})$, and so $F_5\in (W^*(B))^{(2)}$ from Theorem 2.5.

On the contrary, suppose $F=F_1\oplus F_2\oplus F_3\oplus F_4\oplus F_5$ satisfies:

$(1)$ $F_2\oplus F_{3}=0$ or $I$, and $F_5\in (W^*(B))^{(2)}$ if $\lambda=1$;

$(2)$ $F_5\in (W^*(B))^{(2)}$ if $\lambda\in\Lambda_{1}$.\\
It is apparent that $F\in (W^*(\mathcal{T}_\lambda))^{\prime}$, and so we have
\begin{center}
$\left\{\begin{array}{l} W^*(\mathcal{T}_\lambda)=\mathbb C I\oplus \mathcal{B(N}(T-I))\oplus \mathbb C I\oplus M_{2}(\{B\}^{\prime})\quad\hbox{ if }\lambda=1,\\
W^*(\mathcal{T}_\lambda)=\mathbb C I\oplus \mathbb C I\oplus \mathbb C I\oplus \mathbb C I\oplus M_{2}(\{B\}^{\prime})\quad\hbox{ if }\lambda\in\Lambda_{1}.
\end{array}\right.$
\end{center}
By the double commutant theorem, we get
\begin{center}
$\left\{\begin{array}{l} (W^*(\mathcal{T}_\lambda))^{\prime}=\mathcal{B(N}(T))\oplus\mathbb C I \oplus \mathcal{B(N}(T-2I))\oplus {W^*(B)}^{(2)}\hbox{ if }\lambda=1,\\
(W^*(\mathcal{T}_\lambda))^{\prime}=\mathcal{B(N}(T))\oplus \mathcal{B(N}(T-I))\oplus \mathcal{B(N}(T-\lambda I))\\
\quad\quad\quad\quad\quad\quad\quad~~\oplus \mathcal{B(N}(T-(1+\lambda)I))\oplus {W^*(B)}^{(2)}\hbox{ if }\lambda\in\Lambda_{1}.
\end{array}\right.$
\end{center}\qed\\
{\bf Remark 2.9}\quad{\it  Let $T\in\mathcal{ B(H)}$ be a self-adjoint operator having a generic part $T_0$ with respect to $\lambda\in\Lambda$. If $T$ is the pencil of a pair of projections at $\lambda$, then $W^*(\mathcal{T}_\lambda)$ and $(W^*(\mathcal{T}_\lambda))^{\prime}$ are unitarily equivalent to
$$\left\{\begin{array}{l}\mathbb C I\oplus\mathcal{ B(N}(T-I))\oplus \mathbb C I\oplus M_{2}(\{A_{+}\}^{\prime})\quad\hbox{ if }\lambda=1,\\
\mathbb C I\oplus \mathbb C I\oplus \mathbb C I\oplus \mathbb C I\oplus M_{2}(\{A_{+}\}^{\prime})\quad\hbox{ if }\lambda\in\Lambda_{1}
\end{array}\right.$$
and
$$\left\{\begin{array}{l} \mathcal{B(N}(T))\oplus\mathbb C I \oplus \mathcal{B(N}(T-2I))\oplus {W^*(A_{+})}^{(2)}\quad\hbox{ if }\lambda=1,\\
\mathcal{B(N}(T))\oplus\mathcal{ B(N}(T-I))\oplus\mathcal{ B(N}(T-\lambda I))\oplus\mathcal{ B(N}(T-(1+\lambda)I)) \\
\quad\quad\quad\quad\quad\quad\quad\quad \quad\quad\quad\quad\quad\quad\quad\quad\quad\quad~~\oplus {W^*(A_{+})}^{(2)}\quad\hbox{ if }\lambda\in\Lambda_{1},
\end{array}\right.$$
respectively, where $A=T_{0}-\frac{1+\lambda}{2}I$.}\\

According to Theorem 2.6, we know that if $\lambda\in\Lambda_{1}$, then $\mathcal{T}_{\lambda}$ is connected. Next, we consider the connected components of $\mathcal{T}_\lambda$ if $\lambda=1$. Clearly, $\mathcal{T}_\lambda$ is connected if $\mathcal N(T-I)=0$. Generally, the connected components of $\mathcal{T}_1$ depend on that of $\mathcal{P(N}(T-I))$. From the fact that two projections in $\mathcal{P(H)}$ are connected if and only if they are unitarily equivalent, we have the following result.\\\\
{\bf Corollary 2.10}\quad{\it Let $T\in\mathcal B(\mathcal H)$ be a self-adjoint operator such that $T_{0}-I$ has the form $(2.1)$, where $T_0$ is a generic part of $T$ with respect to $1$.
If $\dim {\mathcal{N}}(T-I)<\infty$, then the connected  components of $\mathcal{T}_1$ are
\begin{center}
${\mathrm{E}}_{k}=\{(P_{1,U}, Q_{1,U})\in\mathcal{ T}_1: \dim {\mathcal{R}}(E)=k\}, \  where \ 0 \leq k \leq  \dim {\mathcal{N}}(T-I).$
\end{center}
If $\dim {\mathcal{N}}(T-I)= \infty$, then the connected  components of $\mathcal{T}_1$ are
\begin{center}
${\mathrm{E}}_{m}=\{(P_{1,U}, Q_{1,U})\in\mathcal{T}_1: \dim {\mathcal{R}}(E)=m\},$
\end{center}
\begin{center}
${\mathrm{(I-E)}}_{n}=\{(P_{1,U}, Q_{1,U})\in\mathcal{T}_1: \dim {\mathcal{R}}(I-E)=n\}$
\end{center} and
\begin{center}
${\mathrm{E}}_{\infty}=\{(P_{1,U}, Q_{1,U})\in\mathcal{T}_1: \dim {\mathcal{R}}(E)=\dim {\mathcal{R}}(I-E)=\infty\},$
\end{center}
where m and n are two arbitrary nonnegative  integers.}\\

\section{The uniqueness of $\lambda$ such that $T=\lambda P+Q$}
Theorem 2.6 gives the general representation of pair $(P, Q)$ of projections such that $T=\lambda P+Q$ for some point $\lambda\in\mathbb{R}$. Using the representation, for self-adjoint operator $T$, we investigate the uniqueness of the real number at which $T$ is the pencils of some pairs of projections.\\\\
{\bf Proposition 3.1}\quad{\it Let $T\in\mathcal{ B(H)}$ be  a self-adjoint operator which is the pencils of pairs of projections at $\lambda, \mu\in\mathbb{R}$, respectively. If one of $\lambda$ and $\mu$ is zero, then another is in $\{-1, 0, 1\}$.}\\\\
{\bf Proof.} Suppose $T=\lambda P+Q=\mu P_{1}+Q_{1}$, where $(P, Q)$ and $(P_{1}, Q_{1})$ are pairs of projections. There is no loss of generality in assuming $\mu=0$. We thus have $T=\lambda P+Q=Q_{1}$, which implies that
$\lambda^{2}P+\lambda PQ+\lambda QP=\lambda P.$

If $\lambda=0$, then $Q=T=Q_{1}$. Otherwise, $PQ+QP=(1-\lambda)P.$ It is easily seen that $PQP+QP=(1-\lambda)P=PQ+PQP,$ which yields $PQ=QP$. From this, it induces that $2PQ=(1-\lambda)P$. If $PQ=0$, there would be $\lambda=1$. If $PQ\neq0$, we would have $\lambda=-1$ from $2PQ=(1-\lambda)PQ$.\qed\\

Assume that self-adjoint operator $T\in\mathcal{ B(H)}$ is the pencils of pairs of projections at $\lambda, \mu\in\mathbb{R}\backslash\{0\}$, respectively. Then there exist pair $(P, Q)$ and pair $(P_{1}, Q_{1})$ of projections such that $T=\lambda P+Q=\mu P_{1}+Q_{1}$. To simplify the following statements, we set $\mathcal{H}_{1}=\mathcal{H}\ominus(\mathcal{N}(T)\oplus \mathcal{N}(T-I))$ and $T_{1}=T|_{\mathcal{H}_{1}}$. If $T_{1}=0$, there would be $T=0\oplus I$. It follows that one of $P=0$ and $Q=0$ holds, which contradicts $P\neq0\neq Q$. Therefore $T_{1}\neq 0$. In combination with Theorem 2.6, we have
$$\left\{
  \begin{array}{ll}
    T_{1}=\lambda I\oplus(1+\lambda)I\oplus (A+\frac{1+\lambda}{2}I)\oplus (-A+\frac{1+\lambda}{2}I),\\
T_{1}= \mu I\oplus(1+\mu)I\oplus (B+\frac{1+\mu}{2}I)\oplus (-B+\frac{1+\mu}{2}I)
  \end{array}
\right.\eqno{(3.1)}$$
in terms of
$$\left\{
  \begin{array}{ll}
    \mathcal{H}_{1}=\mathcal{N}(T-\lambda I)\oplus\mathcal{ N}(T-(1+\lambda)I)\oplus  \mathcal{K}_{1}\oplus \mathcal{K}_{1},\\
\mathcal{H}_{1}=\mathcal{N}(T-\mu I)\oplus \mathcal{N}(T-(1+\mu)I)\oplus \mathcal{K}_{2}\oplus \mathcal{K}_{2},
  \end{array}
\right.\eqno{(3.2)}$$
where
$$\frac{|1-|\lambda||}{2}I<A<\frac{1+|\lambda|}{2}I\hbox{ and }\frac{|1-|\mu||}{2}I<B<\frac{1+|\mu|}{2}I.\eqno{(3.3)}$$\\
{\bf Proposition 3.2}\quad{\it Let $T\in\mathcal{B(H)}$ be a self-adjoint operator. Then there are $\lambda, \mu\in\mathbb{R}\backslash\{0\}$ with $|\lambda-\mu|=1$ such that $T$ is the pencils
of pairs of projections at $\lambda, \mu$ if and only if for some $z\in\mathbb{R}\backslash\{0, 1\}$,
\begin{center}
$\{1, z\}\subseteq\sigma(T)\subseteq\{0, 1, z\}.$
\end{center}
In this case, there are only two real numbers at which $T$ is the pencils of pairs of projections.}\\\\
{\bf Proof.} Assume that $T=\lambda P+Q=\mu P_{1}+Q_{1}$, where $\lambda, \mu\in\mathbb{R}\backslash\{0\}$, $(P, Q)$ and $(P_{1}, Q_{1})$ are pairs of projections. Without loss of generality, suppose $(3.1)$ holds in terms of $(3.2)$ and $\lambda>\mu$, that is $\lambda=1+\mu$ and $\mu\neq-1$.

It is easy to divide into five cases $\mu\in[1, +\infty)$, $\mu\in(0,1)$, $\mu\in(-1,0)$, $\mu\in(-2,-1)$ and $\mu\in(-\infty, -2]$ to consider. Next, we discuss the case $\mu\in[1, +\infty)$, and other cases are similar. Then
\begin{center}
$\frac{\lambda-1}{2}I<A<\frac{1+\lambda}{2}I \hbox{ and } \frac{\mu-1}{2}I<B<\frac{1+\mu}{2}I$
\end{center} by $(3.3)$. This easily forces
$$\left\{
  \begin{array}{ll}
  I\leq(1+\mu)I<A+\frac{1+\lambda}{2}I<(2+\mu)I,~0<-A+\frac{1+\lambda}{2}I<I,\\
 I\leq\mu I<B+\frac{1+\mu}{2}I<(1+\mu)I,~0<-B+\frac{1+\mu}{2}I<I.
  \end{array}
\right.$$
It is immediate that $\mathcal{K}_{1}=\{0\}=\mathcal{K}_{2}$ and $\mathcal{N}(T-(1+\lambda) I)=\{0\}=\mathcal{N}(T-\mu I)$. $(1.2)$ becomes
$T_{1}=(1+\mu)I$ since $T_{1}\neq0$. It follows from the fact that $P, Q, P_{1}, Q_{1}$ aren't zero that
$\{1, \lambda\}\subseteq\sigma(T)\subseteq\{0, 1, \lambda\}$.

 On the contrary, let $\{1, z\}\subseteq\sigma(T)\subseteq\{0, 1, z\}$ for some $z\in\mathbb{R}\backslash\{0, 1\}$. Without loss of generality, we assume
$\sigma(T)=\{0, 1, z\}$. It follows that
$T=0\oplus I\oplus zI$ in terms of $\mathcal{H}=\mathcal{N}(T)\oplus\mathcal{N}(T-I)\oplus\mathcal{N}(T-zI)$. Put $P=0\oplus 0\oplus I, Q=0\oplus I\oplus 0, P_{1}=P=0\oplus 0\oplus I$ and $Q_{1}=0\oplus I\oplus I.$ We easily see
$T=zP+Q=(z-1)P_{1}+Q_{1}$. In fact, it is not hard to find that there are only $z$ and $z-1$ at which $T$ is the pencils of pairs of projections.\qed\\\\
{\bf Proposition 3.3}\quad{\it Let $T\in\mathcal{B(H)}$ be a self-adjoint operator. Then $T$ is the pencils of pairs of projections at non-zero real numbers $\lambda, \mu$ with $\max\{\lambda, \mu\}$=$2\min\{\lambda, \mu\}\neq2$$(resp. \min\{\lambda, \mu\}$=$2\\
\max\{\lambda, \mu\}\neq-2)$ if and only if
\begin{center}$\{1+z, z\}\subseteq\sigma(T)\subseteq\{0, 1, 1+z, z\}$
\end{center}
for some $z\in(0, 1)($$resp$.$z\in(-1, 0))$, $\dim\mathcal{N}(T-(1+z)I)=\dim\mathcal{N}(T-zI).$

At this time, there are only two points at which $T$ is the pencils of pairs of projections.}\\\\
{\bf Proof.}\quad Assume that $T=\lambda P+Q=\mu P_{1}+Q_{1}$, where $\lambda, \mu\in\mathbb{R}\backslash\{0\}$, $(P, Q)$ and $(P_{1}, Q_{1})$ are pairs of projections. Without loss of generality, suppose $(3.1)$ holds in terms of $(3.2)$ and $\lambda=2\mu$. Clearly, $\mu\neq\pm1$.

Firstly, we conclude $\mu\notin(-\infty, -1)\cup(1, +\infty)$. Indeed, if $\mu\in(-\infty, -1)\cup(1, +\infty)$, it remains to consider the case $\mu\in(1, +\infty)$, and the case $\mu\in(-\infty, -1)$ is similar. Then
\begin{center}
$\frac{\lambda-1}{2}I<A<\frac{1+\lambda}{2}I \hbox{ and } \frac{\mu-1}{2}I<B<\frac{1+\mu}{2}I$
\end{center}
by $(3.3)$. By direct computation, we have
$$\left\{
  \begin{array}{ll}
2\mu I<\lambda I<A+\frac{1+\lambda}{2}I<(1+\lambda)I,~0<-A+\frac{1+\lambda}{2}I<I,\\
\mu I<B+\frac{1+\mu}{2}I<(1+\mu)I<2\mu I,~0<-B+\frac{1+\mu}{2}I<I.
 \end{array}
\right.$$
Clearly, $\mathcal{K}_{1}=\{0\}=\mathcal{K}_{2}$. Moreover,
\begin{center}
$\mathcal{N}(T-(1+\lambda)I)\oplus\mathcal{N}(T-\lambda I)=\{0\}=\mathcal{N}(T-\mu I)\oplus\mathcal{N}(T-(1+\mu)I),$
\end{center}this induces that $T_{1}=0$, a contradiction. This establishes the conclusion above.
Next, we naturally divide into four cases $\mu\in[\frac{1}{2},1)$, $\mu\in(0, \frac{1}{2})$, $\mu\in(-\frac{1}{2},0)$ and $\mu\in(-1, \frac{1}{2}]$ to consider. By the symmetry, it is sufficient to discuss the cases $\mu\in[\frac{1}{2},1)$ and $\mu\in(0, \frac{1}{2})$, the cases $\mu\in(-\frac{1}{2},0)$ and $\mu\in(-1, \frac{1}{2}]$ are similar.

Case 1: $\mu\in[\frac{1}{2},1)$. Then $\lambda\in[1,2)$. This makes
\begin{center}
$\frac{\lambda-1}{2}I<A<\frac{1+\lambda}{2}I\hbox{ and } \frac{1-\mu}{2}I<B<\frac{1+\mu}{2}I$
\end{center}
by $(3.3)$. By direct computation, we have
$$\left\{
  \begin{array}{ll}
\lambda I<A+\frac{1+\lambda}{2}I<(1+\lambda)I,~0<-A+\frac{1+\lambda}{2}I<I,\\
I<B+\frac{1+\mu}{2}I<(1+\mu)I,~0<-B+\frac{1+\mu}{2}I<\mu I.
  \end{array}
\right.\eqno{(3.4)}$$
Clearly, $\lambda\in\sigma_{p}(B+\frac{1+\mu}{2}I)$, $1+\mu\in\sigma_{p}(A+\frac{1+\lambda}{2}I)$, $\mu\in\sigma_{p}(-A+\frac{1+\lambda}{2}I)$ and $\mathcal{N}(T-(1+\lambda)I)=\{0\}$ by $(3.4)$. Now $(3.1)$ becomes
$$\left\{
  \begin{array}{ll}
   T_{1}=\lambda I\oplus (A+\frac{1+\lambda}{2}I)\oplus (-A+\frac{1+\lambda}{2}I),\\
T_{1}=\mu I\oplus(1+\mu)I\oplus (B+\frac{1+\mu}{2}I)\oplus (-B+\frac{1+\mu}{2}I),
  \end{array}
\right.\eqno{(3.5)}$$
which easily forces that
\begin{center}
$\left\{
  \begin{array}{ll}
    T_{1}-I=(\lambda-1)I\oplus (A+\frac{\lambda-1}{2}I)\oplus (-A+\frac{\lambda-1}{2}I),\\
T_{1}-I=(\mu-1) I\oplus\mu I\oplus (B+\frac{\mu-1}{2}I)\oplus (-B+\frac{\mu-1}{2}I),
\end{array}
\right.$
\end{center}
where
$$\left\{
  \begin{array}{ll}
  0\leq(\lambda-1)I<A+\frac{\lambda-1}{2}I<\lambda I,~-I<-A+\frac{\lambda-1}{2}I<0,\\
0<B+\frac{\mu-1}{2}I<\mu I<I,~-I<-B+\frac{\mu-1}{2}I<(\mu-1)I<0.
\end{array}
\right.$$ It is not hard to see
$$\left\{
  \begin{array}{ll}
(T_{1}-I)_{+}&=(\lambda-1)I\oplus (A+\frac{\lambda-1}{2}I)\oplus 0\\
&=0\oplus \mu I\oplus(B+\frac{\mu-1}{2}I)\oplus 0,\\
(T_{1}-I)_{-}&=0\oplus 0\oplus (A-\frac{\lambda-1}{2}I)\\
&=(1-\mu)I\oplus0\oplus 0\oplus (B-\frac{\mu-1}{2}I).
\end{array}
\right.\eqno{(3.6)}$$ As $A\pm\frac{\lambda-1}{2}I$ and $B\pm\frac{\mu-1}{2}I$ are dense range,
 we see that $$\mathcal{N}(T-\lambda I)\oplus\mathcal{K}_{1}=\mathcal{N}(T-(1+\mu)I)\oplus\mathcal{K}_{2}\hbox{ and } \mathcal{K}_{1}=\mathcal{N}(T-\mu I)\oplus \mathcal{K}_{2}.\eqno{(3.7)}$$
Since $T_{1}\neq0$, it follows that $\mathcal{K}_{1}\neq\{0\}$ from the fact that $\mu\neq1$.

If $\mu=\frac{1}{2}$, we would know $\mathcal{K}_{2}=\{0\}$. Indeed, if $\mathcal{K}_{2}\neq\{0\}$, then we have $(A+\frac{\lambda-1}{2}I)|_{\mathcal{K}_{2}}=B+\frac{\mu-1}{2}I$ and $(A-\frac{\lambda-1}{2}I)|_{\mathcal{K}_{2}}=B-\frac{\mu-1}{2}I$ according to $(3.6)$ and $(3.7)$, and so $\lambda=\mu$. It contradicts $\mu\neq0$ by $\lambda=2\mu$. Therefore $T_{1}=\frac{1}{2}I\oplus\frac{3}{2}I$ in terms of $\mathcal{H}=\mathcal{K}_{1}\oplus\mathcal{K}_{1}=\mathcal{N}(T-\mu I)\oplus\mathcal{N}(T-(1+\mu)I)$ from $(3.5)$ and $(3.7)$. Moreover, we have $\mathcal{N}(T-\mu I)=\mathcal{N}(T-(1+\mu)I)=\mathcal{K}_{1}$ and
\begin{center}
$\{\mu, 1+\mu\}\subseteq\sigma(T)\subseteq\{\mu, 1+\mu, 0, 1\}.$
\end{center}

Next, we consider the case $\mu\in(\frac{1}{2}, 1)$. Firstly, we claim that $\mathcal{K}_{2}=\{0\}$. Indeed, if $\mathcal{K}_{2}\neq\{0\}$, then $(A-\frac{\lambda-1}{2}I)|_{\mathcal{K}_{2}}=B-\frac{\mu-1}{2}I$. Moreover, one obtains $A+\frac{\lambda-1}{2}I\neq\mu I$ and $B+\frac{\mu-1}{2}I\neq(\lambda-1)I$. Otherwise, there would be $\lambda=\mu+\frac{1}{2}$ from the fact that $(A-\frac{\lambda-1}{2}I)|_{\mathcal{K}_{2}}=B-\frac{\mu-1}{2}I$. This contradicts that $\mu\neq\frac{1}{2}$ by $\lambda=2\mu$. Therefore there is a reducing subspace $\mathcal{M}\subseteq \mathcal{K}_{1}\cap \mathcal{K}_{2}$ of $T$ such that $(A+\frac{\lambda-1}{2}I)|_{\mathcal{M}}=(B+\frac{\mu-1}{2}I)|_{\mathcal{M}}$ combining $(3.6)$ with $(3.7)$. This implies that $\lambda=\mu$, which contradicts $\mu\neq0$ from $\lambda=2\mu$. Since $T_{1}\neq0$, it follows that $\mathcal{N}(T-(1+\mu)I)\neq\{0\}$ and $\mathcal{N}(T-\mu I)\neq\{0\}$ from the claim above and $(3.7)$. One obtains $A-\frac{\lambda-1}{2}I=(1-\mu)I$ and $A+\frac{\lambda-1}{2}I=\mu I$ from $(3.6)$ and $(3.7)$. Moreover, if $\mathcal{N}(T-\lambda I)\neq\{0\}$, there would be $\mu=1$, a contradiction. Therefore $\mathcal{N}(T-\lambda I)=\{0\}$. It induces that $T_{1}=(1+\mu)I\oplus\mu I$ in terms of
$$\mathcal{H}=\mathcal{K}_{1}\oplus\mathcal{K}_{1}=\mathcal{N}(T-\mu I)\oplus\mathcal{N}(T-(1+\mu)I)$$
with $\mathcal{K}_{1}=\mathcal{N}(T-\mu I)=\mathcal{N}(T-(1+\mu)I)$ from $(3.7)$.

Case 2: $\mu\in(0, \frac{1}{2})$. Then $\lambda\in(0, 1)$. It makes that $\frac{1-\lambda}{2}I<A<\frac{1+\lambda}{2}I$ and $\frac{1-\mu}{2}I<B<\frac{1+\mu}{2}I$ by $(3.3)$. We easily get
$$\left\{
  \begin{array}{ll}
I<A+\frac{1+\lambda}{2}I<(1+\lambda)I,~0<-A+\frac{1+\lambda}{2}I<\lambda I,\\
I<B+\frac{1+\mu}{2}I<(1+\mu)I,~0<-B+\frac{1+\mu}{2}I<\mu I,
  \end{array}
\right.$$
and consequently $\mathcal{N}(T-\lambda I)\oplus \mathcal{N}(T-(1+\lambda)I)=\{0\}$, $1+\mu\in\sigma_{p}(A+\frac{1+\lambda}{2}I)$ and $\mu\in\sigma_{p}(-A+\frac{1+\lambda}{2}I)$.  Now $(3.1)$ becomes
$$\left\{
  \begin{array}{ll}
    T_{1}=\lambda I\oplus (1+\lambda)I\oplus(A+\frac{1+\lambda}{2}I)\oplus (-A+\frac{1+\lambda}{2}I),\\
T_{1}=\mu I\oplus (1+\mu)I\oplus(B+\frac{1+\mu}{2}I)\oplus (-B+\frac{1+\mu}{2}I)
  \end{array}
\right.\eqno{(3.8)}$$
with respect to $(3.2)$. It follows that

\quad\quad\quad\quad\quad$T_{1}-I=(A-\frac{1-\lambda}{2}I)\oplus (-A-\frac{1-\lambda}{2}I)$

\quad\quad\quad\quad\quad\quad\quad\quad$=(\mu-1)I\oplus \mu I\oplus(B-\frac{1-\mu}{2}I)\oplus (-B-\frac{1-\mu}{2}I),$\\
where
$$\left\{
  \begin{array}{ll}
0<A-\frac{1-\lambda}{2}I<\lambda I,~-I<-A-\frac{1-\lambda}{2}I<(\lambda-1)I<0,\\
 0<B-\frac{1-\mu}{2}I<\mu I,~-I<-B-\frac{1-\mu}{2}I<(\mu-1)I<0.
  \end{array}
\right.$$
Obviously,
$$\left\{
  \begin{array}{ll}
  (T_{1}-I)_{+}=(A-\frac{1-\lambda}{2}I)\oplus 0=0\oplus \mu I\oplus (B-\frac{1-\mu}{2}I)\oplus0  ,\\
 (T_{1}-I)_{-}=0\oplus (A+\frac{1-\lambda}{2}I)=(1-\mu)I\oplus0\oplus 0\oplus (B+\frac{1-\mu}{2}I).
  \end{array}
\right.\eqno{(3.9)}$$
Note that $A+\frac{1-\lambda}{2}I\neq A-\frac{1-\lambda}{2}I$ since $\lambda=2\mu<1$. From the fact that $A\pm\frac{1-\lambda}{2}I$ and $B\pm\frac{1-\mu}{2}I$ are dense range, we easily see that
$$\mathcal{K}_{1}=\mathcal{N}(T-(1+\mu)I)\oplus \mathcal{K}_{2}\hbox{ and } \mathcal{K}_{1}=\mathcal{N}(T-\mu I)\oplus \mathcal{K}_{2} .\eqno{(3.10)}$$
Since $T_{1}\neq0$, we see that $\mathcal{K}_{1}\neq\{0\}$. In fact, $\mathcal{K}_{2}\neq\{0\}$. Otherwise, there is a reducing subspace $\mathcal{M}\subseteq \mathcal{K}_{2}$ of $T$ such that $(B-\frac{1-\mu}{2}I)|_{\mathcal{M}}=(A-\frac{1-\lambda}{2}I)|_{\mathcal{M}}$ and $(A+\frac{1-\lambda}{2}I)|_{\mathcal{M}}=(B+\frac{1-\mu}{2}I)|_{\mathcal{M}}$
from $(3.9)$ and $(3.10)$, which induces that $(B+\frac{\mu-\lambda}{2}I)|_{\mathcal{M}}=A|_{\mathcal{M}}=(B-\frac{\mu-\lambda}{2}I)|_{\mathcal{M}}$. Therefore $\lambda=\mu$. It contradicts $\mu\neq0$ by $\lambda=2\mu$. One may have $A-\frac{1-\lambda}{2}I=\mu I$ and $A+\frac{1-\lambda}{2}I=(1-\mu)I$ from $(3.9)$ and $(3.10)$. It is evident that $T_{1}=\mu I\oplus(1+\mu)I$
in terms of $$\mathcal{H}_{1}=\mathcal{K}_{1}\oplus\mathcal{K}_{1}=\mathcal{N}(T-\mu I)=\mathcal{N}(T-(1+\mu)I)$$
satisfying $\mathcal{K}_{1}=\mathcal{N}(T-\mu I)=\mathcal{N}(T-(1+\mu)I)$ from $(3.8)$ and $(3.10)$.

Therefore, we know that there is $z\in(0, 1)($$resp$. $z\in(-1, 0))$ such that $$\{1+z, z\}\subseteq\sigma(T)\subseteq\{0, 1, 1+z, z\}$$
 with $\dim\mathcal{N}(T-(1+z)I)=\dim\mathcal{N}(T-zI).$

On the contrary, if $\{1+z, z\}\subseteq\sigma(T)\subseteq\{0, 1, 1+z, z\}$ for $z\in(0, 1)$ and $\dim\mathcal{N}(T-(1+z)I)=\dim\mathcal{N}(T-zI)$, we would know that $\mathcal{N}(T-(1+z)I)$ and $\mathcal{N}(T-zI)$ are unitarily equivalent. Without loss of generality, suppose $\sigma(T)=\{0, 1, 1+z, z\}$ and $\mathcal{K}=\mathcal{N}(T-zI)=\mathcal{N}(T-(1+z)I)$. Obviously, we have $T=0\oplus I\oplus(1+z)I\oplus zI$ in terms of $\mathcal{H}=\mathcal{N}(T)\oplus\mathcal{N}(T-I)\oplus\mathcal{K}\oplus\mathcal{K}$.

Suppose $T=\lambda P+Q$ for some real number $\lambda$ and some pair $(P, Q)$ of projections.
It is easily seen that
\begin{center}
$T=\lambda P+Q=0\oplus I\oplus(1+\lambda)I\oplus\lambda I\oplus(A+\frac{1+\lambda}{2}I)\oplus(-A+\frac{1+\lambda}{2}I)$
\end{center}
with respect to $\mathcal{H}=\mathcal{N}(T)\oplus\mathcal{N}(T-I)\oplus\mathcal{N}(T-(1+\lambda)I)\oplus\mathcal{N}(T-\lambda I)\oplus\mathcal{K}_{1}\oplus\mathcal{K}_{1}$ from $(1.1)$, where $\frac{|1-|\lambda||}{2}I<A<\frac{1+|\lambda|}{2}I$. It immediate that
$$\left\{
  \begin{array}{ll}
(1+\lambda)I<\frac{|1-|\lambda||+1+\lambda}{2}I<A+\frac{1+\lambda}{2}I<\frac{2+|\lambda|+\lambda}{2}I,\\
\frac{\lambda-|\lambda|}{2}I<-A+\frac{1+\lambda}{2}I<\frac{-|1-|\lambda||+1+\lambda}{2}I<(1+\lambda)I.
\end{array}
\right.$$
If $\mathcal{K}=\{0\}$, we would have that $$P=0\oplus0\oplus I\oplus I \hbox{ and } Q=0\oplus I\oplus I\oplus0$$
and $T=zP_{1}+Q_{1}$. If $\mathcal{K}\neq\{0\}$, there would be $A+\frac{1+\lambda}{2}I=(1+z)I$ and $-A+\frac{1+\lambda}{2}I=zI$ from the fact that $-A+\frac{1+\lambda}{2}I\leq A+\frac{1+\lambda}{2}I$.
It induces that $\lambda=2\mu$ and $A=\frac{1}{2}I$. According to $(2.3)$, it is easy to get that
$$P=0\oplus 0\oplus\left(\begin{array}{cccc}\frac{1+z}{2}&\frac{\sqrt{1-z^{2}}}{2}\\
\frac{\sqrt{1-z^{2}}}{2}&\frac{1-z}{2}\\\end{array}\right), Q=0\oplus I\oplus\left(\begin{array}{cccc}1-z^{2}&-z\sqrt{1-z^{2}}\\
-z\sqrt{1-z^{2}}&z^{2}\\\end{array}\right).$$
Moreover, we easily check that $(P_{1}, Q_{1})$ and $(P, Q)$ are two pairs of projections. In addition, $\mathcal{R}(P)\cap \mathcal{R}(Q)=\{0\}=\mathcal{R}(P)\cap \mathcal{N}(Q)$ and $(P|_{\mathcal{K}\oplus\mathcal{K}}, Q|_{\mathcal{K}\oplus\mathcal{K}})$ is a generic part of the pair $(P, Q)$.

From what has been discussed, we know that there are only $z$ and $2z$ at which $T$ is the pencils of pairs of projections.

If $\{1+z, z\}\subseteq\sigma(T)\subseteq\{0, 1, 1+z, z\}$ for $z\in(-1, 0)$ and $\dim\mathcal{N}(T-(1+z)I)=\dim\mathcal{N}(T-zI)$, then we also show that there are only $z$ and $2z$ at which $T$ is the pencils of pairs of projections, which completes the proof.\qed\\\\
{\bf Proposition 3.4}\quad{\it Let $T\in\mathcal{B(H)}$ be a self-adjoint operator. Then $T$ is the pencils
of pairs of projections at non-zero real numbers $\lambda, \mu$ with $|\lambda-\mu|=\frac{1}{2}$, $\max\{\lambda, \mu\}\neq1$ and $\min\{\lambda, \mu\}\neq-1$ if and only if one of the followings holds:

$a)$ $\{\frac{1}{2},\frac{1}{2}+z,1+z\}\subseteq\sigma(T)\subseteq\{\frac{1}{2},\frac{1}{2}+z, 1+z, 0, 1\}$ for some $z>\frac{1}{2}$
and $\dim\mathcal{N}(T-\frac{1}{2}I)=\dim\mathcal{N}(T-(\frac{1}{2}+z)I)=\dim\mathcal{N}(T-(1+z)I);$

$b)$ $\{\frac{1}{2}, z, \frac{1}{2}+z\}\subseteq\sigma(T)\subseteq\{\frac{1}{2}, z, \frac{1}{2}+z, 0, 1\}$ for some $z<-\frac{1}{2}$ and
$\dim\mathcal{N}(T-\frac{1}{2}I)=\dim\mathcal{N}(T-(\frac{1}{2}+z)I)=\dim\mathcal{N}(T-zI)$.

In this case, there are only two real numbers at which $T$ is the pencils of pairs of projections.}\\\\
{\bf Proof.}  Assume that $T=\lambda P+Q=\mu P_{1}+Q_{1}$, where $\lambda, \mu\in\mathbb{R}\backslash\{0\}$, $(P, Q)$ and $(P_{1}, Q_{1})$ are pairs of projections. Without loss of generality, suppose $(3.1)$ holds in terms of $(3.2)$ and $\lambda>\mu$, that is $\lambda=\mu+\frac{1}{2}$ and $\mu\neq\pm\frac{1}{2}$.

Firstly, we conclude that $\mu\notin(-1,-\frac{1}{2})\cup(-\frac{1}{2}, 0)\cup(0,\frac{1}{2})$. Indeed, if $\mu\in(0,\frac{1}{2})$, then $\lambda\in(\frac{1}{2}, 1)$ and \begin{center}
$\frac{1-\lambda}{2}I<A<\frac{1+\lambda}{2}I,~\frac{1-\mu}{2}I<B<\frac{1+\mu}{2}I$
\end{center} by $(3.3)$. This easily induces
$$\left\{
  \begin{array}{ll}
I<A+\frac{1+\lambda}{2}I<(1+\lambda)I,~0<-A+\frac{1+\lambda}{2}I<\lambda I,\\
I<B+\frac{1+\mu}{2}I<(1+\mu)I,~0<-B+\frac{1+\mu}{2}I<\mu I.
  \end{array}
\right.$$
Clearly, $\mathcal{N}(T-(1+\lambda)I)\oplus\mathcal{N}(T-\lambda I)=\{0\},~1+\mu\in\sigma_{p}(A+\frac{1+\lambda}{2}I)\hbox{ and } \mu\in\sigma_{p}(-A+\frac{1+\lambda}{2}I)$. Now $(3.1)$ becomes
$$\left\{
  \begin{array}{ll}
    T_{1}=(A+\frac{1+\lambda}{2}I)\oplus (-A+\frac{1+\lambda}{2}I),\\
T_{1}=\mu I\oplus (1+\mu)I\oplus (B+\frac{1+\mu}{2}I)\oplus (-B+\frac{1+\mu}{2}I),
  \end{array}
\right.$$
and hence
$$\left\{
  \begin{array}{ll}
    T_{1}-I=(A+\frac{\lambda-1}{2}I)\oplus (-A+\frac{\lambda-1}{2}I),\\
T_{1}-I=(\mu-1)I\oplus\mu I\oplus (B+\frac{\mu-1}{2}I)\oplus (-B+\frac{\mu-1}{2}I)
  \end{array}
\right.$$
with
$$\left\{
  \begin{array}{ll}
0<A+\frac{\lambda-1}{2}I<\lambda I,~-I<-A+\frac{\lambda-1}{2}I<(\lambda-1)I<0,\\
0<B+\frac{\mu-1}{2}I<\mu I,~-I<-B+\frac{\mu-1}{2}I<(\mu-1)I .
  \end{array}
\right.$$
This yields that
$$\left\{
  \begin{array}{ll}
    (T_{1}-I)_{+}=(A+\frac{\lambda-1}{2}I)\oplus 0=0\oplus\mu I\oplus(B+\frac{\mu-1}{2}I)\oplus 0,\\
 (T_{1}-I)_{-}=0\oplus (A-\frac{\lambda-1}{2}I)=(1-\mu)I\oplus0\oplus (B-\frac{\mu-1}{2}I).
  \end{array}
\right.$$
Since $A\pm\frac{1+\lambda}{2}I$ and $B\pm\frac{1+\mu}{2}I$ are dense range, it follows that
$$\mathcal{K}_{1}=\mathcal{N}(T-(1+\mu)I)\oplus\mathcal{K}_{2} \hbox{ and } \mathcal{K}_{1}=\mathcal{N}(T-\mu I)\oplus \mathcal{K}_{2}.$$
As $T_{1}\neq0$, we have $\mathcal{K}_{1}\neq\{0\}$. Actually, $\mathcal{K}_{2}=\{0\}$. Otherwise, we get that $(A+\frac{\lambda-1}{2}I)|_{\mathcal{K}_{2}}=B+\frac{\mu-1}{2}I$ and $(A-\frac{\lambda-1}{2}I)=\mathcal{K}_{2}=(B-\frac{\mu-1}{2}I)$, which contradicts $\lambda>\mu$. Consequently, $\mu I=A+\frac{\lambda-1}{2}I$ and $(1-\mu)I=A-\frac{\lambda-1}{2}I$. It contradicts $\mu\neq\frac{1}{2}$.

If $\mu\in(-\frac{1}{2}, 0)$, then we easily check that $T_{1}=0$ by calculation, a contradiction. Suppose $\mu\in(-1, -\frac{1}{2})$. Similar to the case $\mu\in(0, \frac{1}{2})$, we have $\mu=-1$ by considering the positive part and negative part of $T_{1}$.

From what has been discussed, we shall divide into four cases to consider in the following.

Case 1:$\mu\in[1, +\infty)$. Then $\lambda\in[\frac{3}{2}, +\infty)$ and
\begin{center}
$\frac{\lambda-1}{2}I<A<\frac{1+\lambda}{2}I,~\frac{\mu-1}{2}I<B<\frac{1+\mu}{2}I$
\end{center} by $(3.3)$.
This easily forces
$$\left\{
  \begin{array}{ll}
\lambda I<A+\frac{1+\lambda}{2}I<(1+\lambda)I,~0<-A+\frac{1+\lambda}{2}I<I,\\
\mu I<B+\frac{1+\mu}{2}I<(1+\mu)I,~0<-B+\frac{1+\mu}{2}I<I.
  \end{array}
\right.$$
It is immediate that $\mathcal{N}(T-(1+\lambda)I)=\{0\}=\mathcal{N}(T-\mu I)$, $1+\mu\in\sigma_{p}(A+\frac{1+\lambda}{2})$ and $\lambda\in\sigma_{p}(B+\frac{1+\mu}{2})$. Now $(3.1)$ becomes
$$\left\{
  \begin{array}{ll}
    T_{1}=\lambda I\oplus (A+\frac{1+\lambda}{2}I)\oplus (-A+\frac{1+\lambda}{2}I),\\
T_{1}=(1+\mu)I\oplus (B+\frac{1+\mu}{2}I)\oplus (-B+\frac{1+\mu}{2}I).
  \end{array}
\right.\eqno{(3.11)}$$
Clearly, $T_{1}-\frac{1+\mu}{2}I=\frac{2\lambda-1-\mu}{2}I\oplus (A+\frac{\lambda-\mu}{2}I)\oplus (-A+\frac{\lambda-\mu}{2}I)= \frac{1+\mu}{2}I\oplus B\oplus(-B),$
where
$0<\frac{2\lambda-1-\mu}{2}I<A+\frac{\lambda-\mu}{2}I<\frac{1+2\lambda-\mu}{2}I,~
-\frac{1+\mu}{2}I<-A+\frac{\lambda-\mu}{2}I<\frac{1-\mu}{2}I.$
It is obvious that
$$\left\{
  \begin{array}{ll}
   (T_{1}-\frac{1+\mu}{2}I)_{+}=\frac{2\lambda-1-\mu}{2}I\oplus (A+\frac{\lambda-\mu}{2}I)\oplus 0= \frac{1+\mu}{2}I\oplus B\oplus 0 ,\\
 (T_{1}-\frac{1+\mu}{2}I)_{-}=0\oplus 0\oplus (A-\frac{\lambda-\mu}{2}I)= 0\oplus 0\oplus B.
  \end{array}
\right.\eqno{(3.12)}$$
Since $A\pm\frac{\lambda-\mu}{2}I$ and $B$ are dense range, we have
$$\mathcal{K}_{1}=\mathcal{K}_{2} \hbox{ and } \mathcal{N}(T-\lambda I)\oplus \mathcal{K}_{1}=\mathcal{N}(T-(1+\mu)I)\oplus \mathcal{K}_{2},\eqno{(3.13)}$$
which means that $\mathcal{N}(T-\lambda I)=\mathcal{N}(T-(1+\mu)I)$. In fact, $\mathcal{K}_{1}\neq\{0\}$. Indeed, if $\mathcal{K}_{1}=\{0\}$, then we have $T_{1}=\lambda I=(1+\mu)I$ from $(3.11)$ and $(3.13)$, and hence $\lambda=1+\mu\geq2$, which contradicts $\lambda=\frac{1}{2}+\mu$.

Next, we claim $A+\frac{\lambda-\mu}{2}I=\frac{1+\mu}{2}I$ and
$B=\frac{2\lambda-1-\mu}{2}I$. Otherwise, we would have that $A-\frac{\lambda-\mu}{2}I=B$ and there exists a reducing subspace $\mathcal{M}\subseteq\mathcal{K}_{1}$ of $T$ such that $(A+\frac{\lambda-\mu}{2}I)|_{\mathcal{M}}=B|_{\mathcal{M}}$ by $(3.12)$ and $(3.13)$. It follows that
\begin{center}
$(A+\frac{\lambda-\mu}{2}I)|_{\mathcal{M}}=B|_{\mathcal{M}}=(A-\frac{\lambda-\mu}{2}I)|_{\mathcal{M}}.$
\end{center}
We easily see that $\lambda=\mu$, which contradicts $\lambda>\mu$. Therefore $\lambda=\mu+\frac{1}{2}$ from the fact that $A-\frac{\lambda-\mu}{2}I=B$. In this case, $(3.11)$ becomes
$$\left\{
  \begin{array}{ll}
    T_{1}=\lambda I\oplus(\frac{1}{2}+\lambda)I\oplus\frac{1}{2}I,\\
T_{1}=(1+\mu)I\oplus(\frac{1}{2}+\mu)I\oplus\frac{1}{2}I
  \end{array}
\right.$$
in terms of $$\left\{
  \begin{array}{ll}
  \mathcal{H}=\mathcal{N}(T-\lambda I)\oplus\mathcal{K}_{1}\oplus\mathcal{K}_{1},\\
  \mathcal{H}=\mathcal{N}(T-(1+\mu)I)\oplus\mathcal{K}_{2}\oplus\mathcal{K}_{2}
  \end{array}
\right.$$
with $\mathcal{N}(T-\lambda I)=\mathcal{N}(T-(1+\mu)I)=\mathcal{K}_{1}=\mathcal{K}_{2}$ and
$$\{1+\mu,\frac{1}{2}+\mu, \frac{1}{2}\}\subseteq\sigma(T)\subseteq\{1+\mu,\frac{1}{2}+\mu, \frac{1}{2}, 0, 1\}.$$

Case 2: $\mu\in(\frac{1}{2},1)$. Then $\lambda\in(1,\frac{3}{2})$ and
\begin{center}
$\frac{\lambda-1}{2}I<A<\frac{1+\lambda}{2}I\hbox{ and } \frac{1-\mu}{2}I<B<\frac{1+\mu}{2}I$
\end{center}
by $(3.3)$. By direct computation, we have
$$\left\{
  \begin{array}{ll}
\lambda I<A+\frac{1+\lambda}{2}I<(1+\lambda)I,~0<-A+\frac{1+\lambda}{2}I<I,\\
I<B+\frac{1+\mu}{2}I<(1+\mu)I,~0<-B+\frac{1+\mu}{2}I<\mu I.
  \end{array}
\right.$$
Clearly, $\mathcal{N}(T-(1+\lambda)I)=\{0\}$, $\lambda\in\sigma_{p}(B+\frac{1+\mu}{2}I)$, $\mu\in\sigma_{p}(-A+\frac{1+\lambda}{2}I)$ and $1+\mu\in\sigma_{p}(A+\frac{1+\lambda}{2}I)$. Note that $(3.1)$ becomes
$$\left\{
  \begin{array}{ll}
   T_{1}=\lambda I\oplus (A+\frac{1+\lambda}{2}I)\oplus (-A+\frac{1+\lambda}{2}I),\\
T_{1}=\mu I\oplus(1+\mu)I\oplus (B+\frac{1+\mu}{2}I)\oplus (-B+\frac{1+\mu}{2}I).
  \end{array}
\right.\eqno{(3.14)}$$
which easily forces that
$$\left\{
  \begin{array}{ll}
    T_{1}-I=(\lambda-1)I\oplus (A+\frac{\lambda-1}{2}I)\oplus (-A+\frac{\lambda-1}{2}I),\\
T_{1}-I=(\mu-1) I\oplus\mu I\oplus (B+\frac{\mu-1}{2}I)\oplus (-B+\frac{\mu-1}{2}I),
\end{array}
\right.$$
where
$$\left\{
  \begin{array}{ll}
  0\leq(\lambda-1)I<A+\frac{\lambda-1}{2}I<\lambda I,~-I<-A+\frac{\lambda-1}{2}I<0,\\
0<B+\frac{\mu-1}{2}I<\mu I,~-I<-B+\frac{\mu-1}{2}I<(\mu-1)I<0.
\end{array}
\right.$$ It is not hard to see
$$\left\{
  \begin{array}{ll}
(T_{1}-I)_{+}&=(\lambda-1)I\oplus (A+\frac{\lambda-1}{2}I)\oplus 0\\
&=0\oplus \mu I\oplus(B+\frac{\mu-1}{2}I)\oplus 0,\\
(T_{1}-I)_{-}&=0\oplus 0\oplus (A-\frac{\lambda-1}{2}I)\\
&=(1-\mu)I\oplus0\oplus 0\oplus (B-\frac{\mu-1}{2}I).
\end{array}
\right.\eqno{(3.15)}$$ Since $A\pm\frac{\lambda-1}{2}I$ and $B\pm\frac{\mu-1}{2}I$ are dense range,
 we see that $$\mathcal{N}(T-\lambda I)\oplus\mathcal{K}_{1}=\mathcal{N}(T-(1+\mu)I)\oplus\mathcal{K}_{2}\hbox{ and } \mathcal{K}_{1}=\mathcal{N}(T-\mu I)\oplus \mathcal{K}_{2}.\eqno{(3.16)}$$
This implies that $\mathcal{K}_{1}\neq\{0\}$ since $T_{1}\neq0$ and $\mu\neq\frac{1}{2}$. In addition, it is easily seen that $\mathcal{K}_{2}\neq\{0\}$ from $\mu\neq\frac{1}{2}$ again.

Firstly, we claim that $\mathcal{N}(T-(1+\mu)I)\neq\{0\}$ and $\mathcal{N}(T-\lambda I)\neq\{0\}$. If one of $\mathcal{N}(T-(1+\mu)I)=\{0\}$ and $\mathcal{N}(T-\lambda I)=\{0\}$ holds, then there is a reducing subspace $\mathcal{M}\subseteq \mathcal{K}_{1}\cap \mathcal{K}_{2}$ of $T$ such that $(A+\frac{\lambda-1}{2}I)|_{\mathcal{M}}=(B+\frac{\mu-1}{2}I)|_{\mathcal{M}}$ and $(A-\frac{\lambda-1}{2}I)|_{M}=(B-\frac{\mu-1}{2}I)|_{M}$  combining $(3.15)$ with $(3.16)$, which contradicts $\lambda>\mu$.
Therefore $A+\frac{\lambda-1}{2}I=\mu I$ and $B+\frac{\mu-1}{2}I=(\lambda-1)I$. In this case, we easily check $\mathcal{N}(T-\mu I)=\{0\}$ by $\mu\neq\frac{1}{2}$. Hence $(3.14)$ becomes
$$\left\{
  \begin{array}{ll}
    T_{1}=\lambda I\oplus(\frac{1}{2}+\lambda)I\oplus\frac{1}{2}I,\\
T_{1}=(1+\mu)I\oplus(\frac{1}{2}+\mu)I\oplus\frac{1}{2}I
  \end{array}
\right.$$
in terms of $$\left\{
  \begin{array}{ll}
  \mathcal{H}=\mathcal{N}(T-\lambda I)\oplus\mathcal{K}_{1}\oplus\mathcal{K}_{1},\\
  \mathcal{H}=\mathcal{N}(T-(1+\mu)I)\oplus\mathcal{K}_{2}\oplus\mathcal{K}_{2}
  \end{array}
\right.$$
with $\mathcal{N}(T-\lambda I)=\mathcal{N}(T-(1+\mu)I)=\mathcal{K}_{1}=\mathcal{K}_{2}$ and
\begin{center}
$\{1+\mu,\frac{1}{2}+\mu, \frac{1}{2}\}\subseteq\sigma(T)\subseteq\{1+\mu,\frac{1}{2}+\mu, \frac{1}{2}, 0, 1\}.$
\end{center}

Case 3: $\mu\in(-\frac{3}{2},-1)$. It yields $\lambda\in(-1,-\frac{1}{2})$, and then we know that $T$ has the form as in Proposition 3.4, which is similar to Case 2 by considering the positive part and negative part of $T_{1}$.

Case 4: $\mu\in(-\infty, -\frac{3}{2}]$. It is easy to see $\lambda\in(-\infty, -1]$. Similar to Case 1, since $\mu\neq-1$, we see that $\lambda=\mu+\frac{1}{2}$ by considering the positive part and negative part of $T_{1}$.

On the contrary, if $(a)$ holds for $z>\frac{1}{2}$, then $\mathcal{N}(T-\frac{1}{2}I)$, $\mathcal{N}(T-(\frac{1}{2}+z)I)$, and $\mathcal{N}(T-(1+z)I)$ are all unitarily equivalent. Set $\mathcal{K}=\mathcal{N}(T-\frac{1}{2}I)$.
Without loss of generality, we may assume $\sigma(T)=\{\frac{1}{2},\frac{1}{2}+z, 1+z, 0, 1\}$. Then
$$T=0\oplus I\oplus\frac{1}{2}I\oplus(\frac{1}{2}+z)I\oplus(1+z)I\eqno{(3.17)}$$
with respect to $$\mathcal{H}=\mathcal{N}(T)\oplus\mathcal{N}(T-I)\oplus\mathcal{K}\oplus\mathcal{K}\oplus\mathcal{K}.\eqno{(3.18)}$$
From Proposition 3.1, we suppose $T=\beta P+Q$ with $\beta\in\mathbb{R}\backslash\{0\}$ for pair $(P, Q)$ of projections from $(3.17)$. It is easily seen that
\begin{center}
$T=\beta P+Q=0\oplus I\oplus(1+\beta)I\oplus\beta I\oplus(A+\frac{1+\beta}{2}I)\oplus(-A+\frac{1+\beta}{2}I)$
\end{center}
with respect to $\mathcal{H}=\mathcal{N}(T)\oplus\mathcal{N}(T-I)\oplus\mathcal{N}(T-(1+\beta)I)\oplus\mathcal{N}(T-\beta I)\oplus\mathcal{K}_{1}\oplus\mathcal{K}_{1}$ from $(1.1)$, where $\frac{|1-|\beta||}{2}I<A<\frac{1+|\beta|}{2}I$. Next, we claim that there are only $\beta=\frac{1}{2}+z$ and $\beta=z$ at which $T$ is the pencils for some pairs of projections. Indeed, it is evident that $\mathcal{K}_{1}\neq\{0\}$ from $(3.17)$. From the fact that $-A+\frac{1+\beta}{2}I\leq A+\frac{1+\beta}{2}I$, one may consider from the following aspects.

$(a.1)$ If $A+\frac{1+\beta}{2}I=(1+z)I$, then $-A+\frac{1+\beta}{2}I=-(1+z)I+(1+\beta)I$.

$(a.1.1)$ If $-A+\frac{1+\beta}{2}I=(\frac{1}{2}+z)I$, there would be $\beta=\frac{1}{2}+2z$ from $-A+\frac{1+\beta}{2}I=-(1+z)I+(1+\beta)I$. Hence $1+\beta=\frac{3}{2}+2z$. It contradicts $\frac{1}{2}\in\sigma_{p}(T)$.

$(a.1.2)$ If $-A+\frac{1+\beta}{2}I=\frac{1}{2}I$, it yields $\beta=\frac{1}{2}+z$. By an elementary calculation, it is not hard to get that
$$P=0\oplus 0\oplus I\oplus\left(\begin{array}{cccc}\frac{4\beta^{2}-1}{4\beta^{2}}&\frac{\sqrt{4\beta^{2}-1}}{4\beta^{2}}\\
\frac{\sqrt{4\beta^{2}-1}}{4\beta^{2}}&\frac{1}{4\beta^{2}}\\\end{array}\right), Q=0\oplus I\oplus0\oplus\left(\begin{array}{cccc}\frac{2\beta+1}{4\beta}&-\frac{\sqrt{4\beta^{2}-1}}{4\beta}\\
-\frac{\sqrt{4\beta^{2}-1}}{4\beta}&\frac{2\beta-1}{4\beta}\\\end{array}\right)$$
in terms of $(3.18)$. We easily check that $(P, Q)$ is a pair of projections and $\mathcal{R}(P)\cap \mathcal{R}(Q)=\{0\}$. Put
$$P_{0}=\left(\begin{array}{cccc}\frac{4\beta^{2}-1}{4\beta^{2}}&\frac{\sqrt{4\beta^{2}-1}}{4\beta^{2}}\\
\frac{\sqrt{4\beta^{2}-1}}{4\beta^{2}}&\frac{1}{4\beta^{2}}\\\end{array}\right)\hbox{ and }Q_{0}=\left(\begin{array}{cccc}\frac{2\beta+1}{4\beta}&-\frac{\sqrt{4\beta^{2}-1}}{4\beta}\\
-\frac{\sqrt{4\beta^{2}-1}}{4\beta}&\frac{2\beta-1}{4\beta}\\\end{array}\right).$$ It is immediate that the pair $(P_{0}, Q_{0})$ is a generic part of the pair $(P, Q)$ with respect to $\frac{1}{2}+z$.

$(a.2)$ If $A+\frac{1+\beta}{2}I=(\frac{1}{2}+z)I$, then $-A+\frac{1+\beta}{2}I=-(\frac{1}{2}+z)I+(1+\beta)I$. Since $-A+\frac{1+\beta}{2}I<A+\frac{1+\beta}{2}I$, it follows that $-A+\frac{1+\beta}{2}I=\frac{1}{2}I$. Hence $\beta=z$. In this case, put
$$P=0\oplus 0\oplus I\oplus\left(\begin{array}{cccc}\frac{4z^{2}-1}{4z^{2}}&\frac{\sqrt{4z^{2}-1}}{4z^{2}}\\
\frac{\sqrt{4z^{2}-1}}{4z^{2}}&\frac{1}{4z^{2}}\\\end{array}\right), Q=0\oplus I\oplus I\oplus\left(\begin{array}{cccc}\frac{2z+1}{4z}&-\frac{\sqrt{4z^{2}-1}}{4z}\\
-\frac{\sqrt{4z^{2}-1}}{4z}&\frac{2z-1}{4z}\\\end{array}\right)$$
in terms of $(3.18)$. It is evident that
$\mathcal{R}(P)\cap \mathcal{N}(Q)=\{0\}$ and $(P, Q)$ is a pair of projections such that $T=zP+Q$. Put
$$P_{0}=\left(\begin{array}{cccc}\frac{4z^{2}-1}{4z^{2}}&\frac{\sqrt{4z^{2}-1}}{4z^{2}}\\
\frac{\sqrt{4z^{2}-1}}{4z^{2}}&\frac{1}{4z^{2}}\\\end{array}\right)\hbox{ and }
Q_{0}=\left(\begin{array}{cccc}\frac{2z+1}{4z}&-\frac{\sqrt{4z^{2}-1}}{4z}\\
-\frac{\sqrt{4z^{2}-1}}{4z}&\frac{2z-1}{4z}\\\end{array}\right).$$
Apparently, the pair $(P_{0}, Q_{0})$ is a generic part of the pair $(P, Q)$ with respect to $z$, which finishes the proof of the claim.

Similar to $(a)$, if $(b)$ holds, we also know that there are only two real numbers $z$ and $z-\frac{1}{2}$ at which $T$ is the pencils for some pairs of projections.\qed\\

For operator $T$ as in Propositions 3.1, 3.2, 3.3 and 3.4, there are only two real numbers at which $T$ is the pencils for some pairs of projections. Based on this, we shall present the main result of this secton below.\\\\
{\bf Theorem 3.5}\quad{\it Let $T\in\mathcal{B(H)}$ be a self-adjoint operator. Then there are at most two real numbers at which $T$ is the pencils of pairs of projections. Moreover, if $T$ doesn't have the forms as in Propositions 3.1, 3.2, 3.3 and 3.4, then there is at most a real number at which $T$ is the pencil for some pair of projections.}\\\\
{\bf Proof.}  Assume that $T=\lambda P+Q=\mu P_{1}+Q_{1}$, where $\lambda, \mu\in\mathbb{R}$, $(P, Q)$ and $(P_{1}, Q_{1})$ are pairs of projections. Without loss of generality, suppose $\lambda\geq\mu$. What is left is to consider $T$ satisfying the case besides Propositions 3.1, 3.2, 3.3 and 3.4. Next, it is easy to be classified into two steps to discuss.

{\bf Step 1:} There aren't real numbers $\lambda, \mu$ with $\lambda\cdot\mu<0$ at which such $T$ is the pencils of pairs of projections. Here are four cases to discuss.

Case 1: $\lambda\in[1,+\infty)$ and $\mu\in(-1,0)$. Then
\begin{center}
$\frac{\lambda-1}{2}I<A<\frac{1+\lambda}{2}I \hbox{ and } \frac{1+\mu}{2}I<B<\frac{1-\mu}{2}I$
\end{center}
by $(3.3)$. This gives that
$$\left\{
  \begin{array}{ll}
   I\leq\lambda I<A+\frac{1+\lambda}{2}I<(1+\lambda)I,\quad 0<-A+\frac{1+\lambda}{2}I<I,\\
0<(1+\mu)I<B+\frac{1+\mu}{2}I<I,\quad~\mu I<-B+\frac{1+\mu}{2}I<0.
  \end{array}
\right.$$%from $(3.3)$,
Clearly, $\mathcal{K}_{1}=\{0\}=\mathcal{K}_{2}$, $\mathcal{N}(T-\mu I)=\{0\}=\mathcal{N}(T-(1+\lambda)I)$. Moreover, $\mathcal{N}(T-\lambda I)=\{0\}$ if $\lambda>1$ and $T_{1}|_{\mathcal{N}(T-\lambda I)}=\{0\}$ if $\lambda=1$. $(3.1)$ becomes $T_{1}=0$, which contradicts $T_{1}\neq0$.

Case 2: $\lambda\in[1,+\infty)$ and $\mu\in(-\infty,-1]$. We also get a contradiction similar to Case 1.

Case 3: $\lambda\in(0,1)$ and $\mu\in(-1,0)$. This clearly forces
\begin{center}
$\frac{1-\lambda}{2}I<A<\frac{1+\lambda}{2}I \hbox{ and } \frac{1+\mu}{2}I<B<\frac{1-\mu}{2}I$
\end{center}
by $(3.3)$. From this, one knows that
$$\left\{
  \begin{array}{ll}
I<A+\frac{1+\lambda}{2}I<(1+\lambda)I,~~~~~~~0<-A+\frac{1+\lambda}{2}I<\lambda I<I,\\
0<(1+\mu)I<B+\frac{1+\mu}{2}I<I,~\mu I<-B+\frac{1+\mu}{2}I<0.
  \end{array}
\right.$$%from $(3.4)$,
Apparently, $\mathcal{K}_{1}=\{0\}=\mathcal{K}_{2}$ and $\mathcal{N}(T-\mu I)=\{0\}$ and $\mathcal{N}(T-(1+\lambda)I)=\{0\}$. Now $(3.1)$ becomes $T_{1}=(1+\mu)I=\lambda I$, which is the form as in Proposition 3.2, a contradiction.

Case 4: $\lambda\in(0,1)$ and $\mu\in(-\infty,-1]$. Similarly, it contradicts $T_{1}\neq0$.

{\bf Step 2:} If there are two real numbers $\lambda, \mu$ with $\lambda\cdot\mu>0$ at which such $T$ is the pencils of pairs of projections, there would be $\lambda=\mu$. Here are six cases to discuss.

Case 5:~$\lambda, \mu\in[1, +\infty)$. Then $\frac{\lambda-1}{2}I<A<\frac{1+\lambda}{2}I \hbox{ and } \frac{\mu-1}{2}I<B<\frac{1+\mu}{2}I$ by $(3.3)$. This easily forces
$$\left\{
  \begin{array}{ll}
  I\leq\lambda I<A+\frac{1+\lambda}{2}I<(1+\lambda)I,~0<-A+\frac{1+\lambda}{2}I<I,\\
I\leq\mu I<B+\frac{1+\mu}{2}I<(1+\mu)I,~0<-B+\frac{1+\mu}{2}I<I.
  \end{array}
\right.\eqno{(3.19)}$$
It is clear that $\mathcal{N}(T-(1+\lambda)I)=\{0\}=\mathcal{N}(T-\mu I)$ if $\lambda>\mu$. $(3.1)$ becomes
$$\left\{
  \begin{array}{ll}
    T_{1}=\lambda I\oplus (A+\frac{1+\lambda}{2}I)\oplus (-A+\frac{1+\lambda}{2}I),\\
T_{1}=(1+\mu)I\oplus (B+\frac{1+\mu}{2}I)\oplus (-B+\frac{1+\mu}{2}I).
  \end{array}
\right.\eqno{(3.20)}$$
Clearly, $T_{1}-\frac{1+\mu}{2}I=\frac{2\lambda-1-\mu}{2}I\oplus (A+\frac{\lambda-\mu}{2}I)\oplus (-A+\frac{\lambda-\mu}{2}I)= \frac{1+\mu}{2}I\oplus B\oplus(-B),$
where
$0<\frac{2\lambda-1-\mu}{2}I<A+\frac{\lambda-\mu}{2}I<\frac{1+2\lambda-\mu}{2}I,~
-\frac{1+\mu}{2}I<-A+\frac{\lambda-\mu}{2}I<\frac{1-\mu}{2}I.$
It is obvious that
$$\left\{
  \begin{array}{ll}
   (T_{1}-\frac{1+\mu}{2}I)_{+}=\frac{2\lambda-1-\mu}{2}I\oplus (A+\frac{\lambda-\mu}{2}I)\oplus 0= \frac{\mu+1}{2}I\oplus B\oplus 0 ,\\
 (T_{1}-\frac{1+\mu}{2}I)_{-}=0\oplus 0\oplus (A-\frac{\lambda-\mu}{2}I)= 0\oplus 0\oplus B.
  \end{array}
\right.\eqno{(3.21)}$$
Since $A\pm\frac{\lambda-\mu}{2}I$ and $B$ are dense range, we see that
$$\mathcal{K}_{1}=\mathcal{K}_{2} \hbox{ and } \mathcal{N}(T-\lambda I)\oplus \mathcal{K}_{1}=\mathcal{N}(T-(1+\mu)I)\oplus \mathcal{K}_{2}.\eqno{(3.22)}$$

Firstly, we claim that $\lambda\leq1+\mu$. Indeed, if $\lambda>1+\mu$, there would be $\mathcal{N}(T-(1+\mu)I)=\{0\}$ and $\mathcal{N}(T-\lambda I)=\{0\}$ by $(3.19)$. We easily see that
$$\left\{
  \begin{array}{ll}
   (T_{1}-\frac{1+\mu}{2}I)_{+}= (A+\frac{\lambda-\mu}{2}I)\oplus 0=B\oplus 0 ,\\
 (T_{1}-\frac{1+\mu}{2}I)_{-}=0\oplus (A-\frac{\lambda-\mu}{2}I)=  0\oplus B
  \end{array}
\right.$$
from ${(3.21)}$. This implies that $A+\frac{\lambda-\mu}{2}I=B=A-\frac{\lambda-\mu}{2}I$. Hence $\lambda=\mu$. It contradicts $\lambda>\mu$. Therefore we have $1+\mu\in\sigma_{p}(A+\frac{1+\lambda}{2}I)$ and $\lambda\in\sigma_{p}(B+\frac{1+\mu}{2}I)$ by $(3.19)$. Next, we divide into two aspects to consider.

(5.a) Suppose $\mathcal{N}(T-\lambda I)\neq\{0\}$ and $\mathcal{N}(T-(1+\mu)I)\neq\{0\}$. In this case, we would have $\mathcal{K}_{1}\neq\{0\}$ and $A+\frac{\lambda-\mu}{2}I\neq\frac{\mu+1}{2}I$ and
$B\neq\frac{2\lambda-1-\mu}{2}I$. Indeed, if $\mathcal{K}_{1}=\{0\}$, it is easily seen that $T_{1}=\lambda I=(1+\mu)I$ with $\lambda=1+\mu$ from $(3.20)$ and $(3.22)$, which makes that $T$ has the form as in Proposition 3.2. If $A+\frac{\lambda-\mu}{2}I=\frac{\mu+1}{2}I$ and
$B=\frac{2\lambda-1-\mu}{2}I$, there would be $A=\frac{2\mu+1-\lambda}{2}I=\frac{3\lambda-1-2\mu}{2}I$ from the fact that $A-\frac{\lambda-\mu}{2}I=B$, which induces that $\lambda=\mu+\frac{1}{2}$. At this time, it is easy to check that $T$ has the form as in Proposition 3.4. From what have been discussed, we would have that $A-\frac{\lambda-\mu}{2}I=B$ and there exists a reducing subspace $\mathcal{M}\subseteq\mathcal{K}_{1}$ of $T$ such that $(A+\frac{\lambda-\mu}{2}I)|_{\mathcal{M}}=B|_{\mathcal{M}}$ according to $(3.21)$ and $(3.22)$. It follows that $(A+\frac{\lambda-\mu}{2}I)|_{\mathcal{M}}=B|_{\mathcal{M}}=(A-\frac{\lambda-\mu}{2}I)|_{\mathcal{M}}.$ It is immediate that $\lambda=\mu$, which contradicts $\lambda>\mu$.

$(5.b)$ Suppose $\mathcal{N}(T-\lambda I)=\{0\}$ or $\mathcal{N}(T-(1+\mu)I)=\{0\}$. We would know that $\mathcal{K}_{1}\neq\{0\}$ as $T_{1}\neq0$, this implies that there is a reducing subspace $\mathcal{M}\subseteq \mathcal{K}_{1}$ of $T$ satisfying $(A+\frac{\lambda-\mu}{2}I)|_{\mathcal{M}}=B|_{\mathcal{M}}=(A-\frac{\lambda-\mu}{2}I)|_{\mathcal{M}}$ from $(3.21)$ and $(3.22)$. Obviously, $\lambda=\mu$, which contradicts $\lambda>\mu$.

From above, one can get $\lambda=\mu$.

Case 6: $\lambda\in[1,+\infty)$ and $\mu\in(0,1)$. Then
\begin{center}
$\frac{\lambda-1}{2}I<A<\frac{1+\lambda}{2}I\hbox{ and } \frac{1-\mu}{2}I<B<\frac{1+\mu}{2}I$
\end{center}
by $(3.3)$. By direct computation, we have
$$\left\{
  \begin{array}{ll}
\lambda I<A+\frac{1+\lambda}{2}I<(1+\lambda)I,~0<-A+\frac{1+\lambda}{2}I<I,\\
I<B+\frac{1+\mu}{2}I<(1+\mu)I,~0<-B+\frac{1+\mu}{2}I<\mu I.
  \end{array}
\right.$$
Clearly, $\mathcal{N}(T-(1+\lambda)I)=\{0\}$ and $\mu\in\sigma_{p}(-A+\frac{1+\lambda}{2}I)$. $(3.1)$ becomes
$$\left\{
  \begin{array}{ll}
   T_{1}=\lambda I\oplus (A+\frac{1+\lambda}{2}I)\oplus (-A+\frac{1+\lambda}{2}I),\\
T_{1}=\mu I\oplus(1+\mu)I\oplus (B+\frac{1+\mu}{2}I)\oplus (-B+\frac{1+\mu}{2}I),
  \end{array}
\right.\eqno{(3.23)}$$
which easily forces that
$$\left\{
  \begin{array}{ll}
    T_{1}-I=(\lambda-1)I\oplus (A+\frac{\lambda-1}{2}I)\oplus (-A+\frac{\lambda-1}{2}I),\\
T_{1}-I=(\mu-1) I\oplus\mu I\oplus (B+\frac{\mu-1}{2}I)\oplus (-B+\frac{\mu-1}{2}I),
\end{array}
\right.$$
where
$$\left\{
  \begin{array}{ll}
  0\leq(\lambda-1)I<A+\frac{\lambda-1}{2}I<\lambda I,~-I<-A+\frac{\lambda-1}{2}I<0,\\
0<B+\frac{\mu-1}{2}I<\mu I<I,~-I<-B+\frac{\mu-1}{2}I<(\mu-1)I<0.
\end{array}
\right.\eqno{(3.24)}$$ It is not hard to see
$$\left\{
  \begin{array}{ll}
(T_{1}-I)_{+}&=(\lambda-1)I\oplus (A+\frac{\lambda-1}{2}I)\oplus 0\\
&=0\oplus \mu I\oplus(B+\frac{\mu-1}{2}I)\oplus 0,\\
(T_{1}-I)_{-}&=0\oplus 0\oplus (A-\frac{\lambda-1}{2}I)\\
&=(1-\mu)I\oplus0\oplus 0\oplus (B-\frac{\mu-1}{2}I).
\end{array}
\right.\eqno{(3.25)}$$ Since $A\pm\frac{\lambda-1}{2}I$ and $B\pm\frac{\mu-1}{2}I$ are dense range,
 we see that $$\mathcal{N}(T-\lambda I)\oplus\mathcal{K}_{1}=\mathcal{N}(T-(1+\mu)I)\oplus\mathcal{K}_{2}\hbox{ and } \mathcal{K}_{1}=\mathcal{N}(T-\mu I)\oplus \mathcal{K}_{2}.\eqno{(3.26)}$$
In this case, we shall prove that there are not such $\lambda, \mu$ at which  $T$ is the pencils of pairs of projections.

Firstly, we conclude $\lambda-1\leq\mu$. Indeed, if $\lambda>1+\mu$, then $\mathcal{N}(T-\lambda I)=\{0\}$ and $\mathcal{N}(T-(1+\mu)I)=\{0\}$ by $(3.24)$. It leads to $\mathcal{N}(T-\mu I)=\{0\}$ and
$$\left\{
  \begin{array}{ll}
(T_{1}-I)_{+}=(A+\frac{\lambda-1}{2}I)\oplus 0=0\oplus(B+\frac{\mu-1}{2}I)\oplus 0,\\
(T_{1}-I)_{-}= 0\oplus (A-\frac{\lambda-1}{2}I)= 0\oplus (B-\frac{\mu-1}{2}I)
\end{array}
\right.$$
by $(3.25)$ and $(3.26)$, thus  $A+\frac{\lambda-1}{2}I=B+\frac{\mu-1}{2}I$ and $A-\frac{\lambda-1}{2}I=B-\frac{\mu-1}{2}I$. We assert
$\lambda=\mu$, which contradicts $\lambda>\mu$. Therefore $\lambda\in\sigma_{p}(B+\frac{1+\mu}{2}I)$ and $1+\mu\in\sigma_{p}(A+\frac{1+\lambda}{2}I)$ by $(3.24)$. Next, we divide into two aspects to consider.

$(6.a)$ Suppose $\mathcal{N}(T-(1+\mu)I)\neq\{0\}$ and $\mathcal{N}(T-\lambda I)\neq\{0\}$.

$(6.a.1)$ If $\mathcal{K}_{1}=\{0\}$, then $\lambda=1+\mu$ from $(3.23)$ and $(3.26)$. We have that $T_{1}=\lambda I=(1+\mu)I$ in terms of $\mathcal{H}_{1}$ with $\mathcal{H}_{1}=\mathcal{N}(T-\lambda I)=\mathcal{N}(T-(1+\mu)I)$, which makes that $T$ has the form as in Proposition 3.2, a contradiction.

$(6.a.2)$ If $\mathcal{K}_{2}=\{0\}$, then we get $A-\frac{\lambda-1}{2}I=(1-\mu)I$ and $A+\frac{\lambda-1}{2}I=\mu I=(\lambda-1)I$ from $(3.23)$ and $(3.26)$. By a simple computation, it is clear that $\lambda=2\mu=\mu+1$, which contradicts $\mu<1$.

$(6.a.3)$ If $\mathcal{K}_{1}\neq\{0\}$ and $\mathcal{K}_{2}\neq\{0\}$, $A+\frac{\lambda-1}{2}I=\mu I$ and $B+\frac{\mu-1}{2}I=(\lambda-1)I$, there would be $A=\frac{2\mu+1-\lambda}{2}I$, $B=\frac{2\lambda-1-\mu}{2}I$ and $(A-\frac{\lambda-1}{2}I)|_{\mathcal{K}_{2}}=B-\frac{\mu-1}{2}I$. It ensures $\lambda=\mu+\frac{1}{2}$ and $\mu\geq\frac{1}{2}$ from $\lambda\geq1$. Obviously,  $$\mathcal{N}(T-\lambda I)=\mathcal{N}(T-(1+\mu)I)=\mathcal{K}_{1}=\mathcal{K}_{2},$$ and then we know that $T$ has the form as in Proposition 3.4. It is a contradiction.

$(6.a.4)$ If $\mathcal{K}_{1}\neq\{0\}$ and $\mathcal{K}_{2}\neq\{0\}$, $A+\frac{\lambda-1}{2}I\neq\mu I$ and $B+\frac{\mu-1}{2}I\neq(\lambda-1)I$, then $(A-\frac{\lambda-1}{2}I)|_{\mathcal{K}_{2}}=B-\frac{\mu-1}{2}I$ and there is a reducing subspace $\mathcal{M}\subseteq \mathcal{K}_{2}$ of $T$ such that $(A+\frac{\lambda-1}{2}I)|_{\mathcal{M}}=(B+\frac{\mu-1}{2}I)|_{\mathcal{M}}$ from $(3.25)$ and $(3.26)$. Clearly, $\lambda=\mu$, it contradicts $\lambda>\mu$.

$(6.b)$ Suppose $\mathcal{N}(T-(1+\mu)I)=\{0\}$ or $\mathcal{N}(T-\lambda I)=\{0\}$. Since $T_{1}\neq0$, it follows that $\mathcal{N}(T-(1+\mu)I)=\{0\}$ and $\mathcal{K}_{2}=\{0\}$ cann't hold at the same time from $(3.26)$.

$(6.b.1)$ If $\mathcal{K}_{2}=\{0\}$ and $\mathcal{N}(T-\lambda I)=\{0\}$, then $A-\frac{\lambda-1}{2}I=(1-\mu)I$ and $A+\frac{\lambda-1}{2}I=\mu I$
from $(3.25)$ and $(3.26)$. This means that $\mathcal{K}_{1}=\mathcal{N}(T-\mu I)=\mathcal{N}(T-(1+\mu)I)$ and $\lambda=2\mu$, which forces that $T$ has the form as in Proposition 3.3, a contradiction.

$(6.b.2)$ If $\mathcal{K}_{2}\neq\{0\}$, we would have that $\mathcal{K}_{1}\neq\{0\}$ and there is a reducing subspace $\mathcal{M}\subseteq \mathcal{K}_{1}\cap \mathcal{K}_{2}$ of $T$ such that $(A+\frac{\lambda-1}{2}I)|_{\mathcal{M}}=(B+\frac{\mu-1}{2}I)|_{\mathcal{M}}$ and $(A-\frac{\lambda-1}{2}I)|_{K_{2}}=B-\frac{\mu-1}{2}I$  combining $(3.25)$ with $(3.26)$. This implies that  $(A-\frac{\lambda-\mu}{2}I)|_{\mathcal{M}}=B|_{\mathcal{M}}=(A+\frac{\lambda-\mu}{2}I)|_{\mathcal{M}}$, hence $\lambda=\mu$. It contradicts $\lambda>\mu$.

Case 7: $\lambda, \mu\in(0, 1)$. Then $\frac{1-\lambda}{2}I<A<\frac{1+\lambda}{2}I$ and $\frac{1-\mu}{2}I<B<\frac{1+\mu}{2}I$ by $(3.3)$. We easily get
$$\left\{
  \begin{array}{ll}
I<A+\frac{1+\lambda}{2}I<(1+\lambda)I,~0<-A+\frac{1+\lambda}{2}I<\lambda I,\\
I<B+\frac{1+\mu}{2}I<(1+\mu)I,~0<-B+\frac{1+\mu}{2}I<\mu I.
  \end{array}
\right.$$
In this way, $\mathcal{N}(T-\lambda I)\oplus \mathcal{N}(T-(1+\lambda)I)=\{0\}$, $1+\mu\in\sigma_{p}(A+\frac{1+\lambda}{2}I)$ and $\mu\in\sigma_{p}(-A+\frac{1+\lambda}{2}I)$ if $\lambda>\mu$.  Now $(3.1)$ becomes
$$\left\{
  \begin{array}{ll}
    T_{1}=\lambda I\oplus (1+\lambda)I\oplus(A+\frac{1+\lambda}{2}I)\oplus (-A+\frac{1+\lambda}{2}I),\\
T_{1}=\mu I\oplus (1+\mu)I\oplus(B+\frac{1+\mu}{2}I)\oplus (-B+\frac{1+\mu}{2}I)
  \end{array}
\right.$$
with respect to $(3.2)$. It follows that

\quad\quad\quad\quad\quad$T_{1}-I=(A-\frac{1-\lambda}{2}I)\oplus (-A-\frac{1-\lambda}{2}I)$

\quad\quad\quad\quad\quad\quad\quad\quad$=(\mu-1)I\oplus \mu I\oplus(B-\frac{1-\mu}{2}I)\oplus (-B-\frac{1-\mu}{2}I),$\\
where
\begin{center}$\left\{
  \begin{array}{ll}
0<A-\frac{1-\lambda}{2}I<\lambda I,~-I<-A-\frac{1-\lambda}{2}I<(\lambda-1)I<0,\\
 0<B-\frac{1-\mu}{2}I<\mu I,~-I<-B-\frac{1-\mu}{2}I<(\mu-1)I<0.
  \end{array}
\right.$
\end{center}
Obviously,
$$\left\{
  \begin{array}{ll}
  (T_{1}-I)_{+}=(A-\frac{1-\lambda}{2}I)\oplus 0=0\oplus \mu I\oplus (B-\frac{1-\mu}{2}I)\oplus0  ,\\
 (T_{1}-I)_{-}=0\oplus (A+\frac{1-\lambda}{2}I)=(1-\mu)I\oplus0\oplus 0\oplus (B+\frac{1-\mu}{2}I).
  \end{array}
\right.\eqno{(3.27)}$$
Note that $A+\frac{1-\lambda}{2}I\neq A-\frac{1-\lambda}{2}I$ since $\lambda<1$. From the fact that $A\pm\frac{1-\lambda}{2}I$ and $B\pm\frac{1-\mu}{2}I$ are dense range, we easily see that
$$\mathcal{K}_{1}=\mathcal{N}(T-\mu I)\oplus \mathcal{K}_{2} \hbox{ and } \mathcal{K}_{1}=\mathcal{N}(T-(1+\mu)I)\oplus \mathcal{K}_{2}.\eqno{(3.28)}$$
Since $T_{1}\neq0$, $\mathcal{K}_{1}\neq\{0\}$. In fact, $\mathcal{K}_{2}\neq\{0\}$. Otherwise, we have $A-\frac{1-\lambda}{2}I=\mu I$ and $A+\frac{1-\lambda}{2}I=(1-\mu)I$ from $(3.27)$ and $(3.28)$. It is evident that $\mathcal{K}_{1}=\mathcal{N}(T-\mu I)=\mathcal{N}(T-(1+\mu)I)$ and $\lambda=2\mu$, which makes that $T$ has the form as in Proposition 3.3. It contradicts the assumption. Therefore there is a reducing subspace $\mathcal{M}\subseteq \mathcal{K}_{2}$ of $T$ such that $(B-\frac{1-\mu}{2}I)|_{\mathcal{M}}=(A-\frac{1-\lambda}{2}I)|_{\mathcal{M}}$ and $(A+\frac{1-\lambda}{2}I)|_{\mathcal{M}}=(B+\frac{1-\mu}{2}I)|_{\mathcal{M}}$
from $(3.27)$ and $(3.28)$. It induces $\lambda=\mu$, which contradicts $\lambda>\mu$.

Consequently, $\lambda=\mu$.

Case 8: $\lambda, \mu\in(-1, 0)$. From the symmetry of Case 7, we also get $\lambda=\mu$.

Case 9: $\lambda\in(-1,0)$ and $\mu\in(-\infty,-1]$. This case is similar to Case 6 by considering the positive part and negative part of $T$.

Case 10: $\lambda, \mu\in(-\infty, -1]$. Similar to Case 5, one can have $\lambda=\mu$ by considering the positive part and negative part of $T_{1}$.

From the Case 1 to Case 10, we know that if self-adjoint operator $T$, except for the cases as in Propositions 3.1, 3.2, 3.3 and 3.4, is the pencils of pairs of projections at two real numbers, then the two real numbers are equal, which completes the proof.\qed

%% The Appendices part is started with the command \appendix;
%% appendix sections are then done as normal sections
%% \appendix

%% \section{}
%% \label{}

%% References
%%
%% Following citation commands can be used in the body text:
%% Usage of \cite is as follows:
%%   \cite{key}          ==>>  [#]
%%   \cite[chap. 2]{key} ==>>  [#, chap. 2]
%%   \citet{key}         ==>>  Author [#]

%% References with bibTeX database:

%\bibliographystyle{model1a-num-names}
%\bibliography{<your-bib-database>}

%% Authors are advised to submit their bibtex database files. They are
%% requested to list a bibtex style file in the manuscript if they do
%% not want to use model1a-num-names.bst.

%% References without bibTeX database:

\end{document}